\UseRawInputEncoding

\documentclass[12pt]{amsart}
\usepackage{amssymb,amscd}  
\usepackage{pb-diagram}   
\usepackage{verbatim}

\usepackage[dvipdfmx]{graphicx,color}  

\usepackage{graphicx}
\usepackage{mathrsfs}

\usepackage{here}   

\usepackage[english]{babel}
\usepackage{amsfonts,amssymb,amsmath,amstext,amsbsy,amsopn,amscd,amsthm,graphicx,euscript,graphicx}
\usepackage{graphics}
\usepackage[all]{xy}

\newtheorem*{Whitney towers}{Theorem~\ref{Whitney towers}}
\newtheorem*{h-towers}{Theorems ~\ref{half} \& \ref{$(n)$-solvable}}

\newtheorem*{surgery curves}{Theorem~\ref{surgery curves}}
\newtheorem*{cg=0}{Theorem~\ref{vanish}}

\newtheorem{thm}{Theorem}[section]
\newtheorem{mth}[thm]{Main Theorem}
\newtheorem{proposition}[thm]{Proposition} 

\newtheorem{corollary}[thm]{Corollary}

\theoremstyle{definition}
\newtheorem{definition}[thm]{Definition}

\newtheorem{remark}[thm]{Remark}

\numberwithin{equation}{section}
\numberwithin{figure}{section}
\numberwithin{table}{section}

\newcommand{\vs}{\vskip10mm}
\newcommand{\bb}{\bigbreak}

\newcommand{\h}{\noindent}
\newcommand{\x}{\times}

\newcommand{\Z}{\mathbb{Z}}

\newcommand{\C}{\mathbb{C}}

\newcommand{\R}{\mathbb{R}}

\def\yen{{\setbox0=\hbox{Y}Y\kern-.97\wd0\vbox{hrule height.lex width.98%
\wd0\kern.33ex\hrule height.lex width.98\wd0\kern.45ex}}}

\begin{document}
{

{\bf 
\h Quantum Invariants of Links and 3-Manifolds with Boundary defined via Virtual Links:
Calculation of some examples
}
\\

\vskip7mm
\h
Heather A. Dye,   
Louis H. Kauffman and Eiji Ogasa

\vskip14mm


\h{\bf Abstract.}  
In the prequel of this paper, 
Kauffman and Ogasa introduced new topological quantum invariants 
of compact oriented 3-manifolds with boundary 
where the boundary is a disjoint union of two identical surfaces. 
The invariants are constructed via surgery on manifolds of the form $F \times I$ where $I$ denotes the unit interval.
Since virtual knots and links are represented as links in such thickened surfaces, we are able also to construct invariants in terms of 
virtual link diagrams (planar diagrams with virtual crossings).\\
 
These invariants are new, nontrivial, and calculable 
examples of quantum invariants of 
3-manifolds with non-vacuous boundary.\\

Since virtual knots and links are represented by embeddings of circles in thickened surfaces, we refer to embeddings of circles in the 3-sphere as {\it classical links}.
Classical links are the same as virtual links that can be represented in a thickened 2-sphere and it is a fact that classical links, up to isotopy, embed in the collection of virtual links taken up
to isotopy.
We give a new invariant of classical links in the 3-sphere in the following sense: 
Consider a link $L$ in $S^3$ of two components. 
The complement of a tubular neighborhood of $L$ 
is a manifold whose boundary consists in two copies of a torus. 
Our invariants apply to this case of bounded manifold 
and give new invariants of the given link of two components. 
Invariants of knots are also obtained. 
\\


In this paper we calculate the topological quantum invariants of  some examples  explicitly.  
We conclude from our examples that our invariant is new and strong enough to distinguish some classical knots from one another.  

We examine links that are embedded in thickened surfaces 
and obtain invariants of three manifolds obtained by surgery 
on these thickened surfaces. 
One could take the viewpoint that the thickened surfaces are 
embedded in the three sphere and so also consider the three manifolds 
obtained by surgery on the links in the three sphere. 
These two points of view are distinct and give distinct invariants. 
This point of view for links in thickened surfaces 
is distinct 
from the usual point of view for the Reshetikhin-Turaev invariants, 
and our invariants give distinct results from these invariants.
(See the body what kind of viewpoint).

\tableofcontents

\section{
\bf Introduction}\label{psecintro}
In the prequel \cite{KOq} of this paper, 
Kauffman and Ogasa introduced new topological quantum invariants 
of compact oriented 3-manifolds with boundary 
where the boundary is a disjoint union of two identical surfaces, 
by using {\it virtual links}. 
Let 
$\upsilon_r$ 
denote this topological invariant. 
Using the topological quantum invariants 
$\upsilon_r$, 
they defined invariants of classical knots and links in the 3-sphere.  
%
We review the invariant 
$\upsilon_r$ 
in Part \ref{part1}. 

The set of virtual links is a quotient set of links in thickened surfaces 
(\cite{Kauffman1,Kauffman,Kauffmani}). 
We review it in \S\ref{psecv1k}.  
In this paper, a thickened surface means 
(an oriented closed surface) $\x$(the oriented interval). \\
Links in the 3-sphere are called {\it classical links}. 
There is a natural bijection between 
the set of links in the 3-spheres and that in the thickened 2-sphere.
Therefore the set of classical links is a subset of that of virtual links (In fact, it is a proper subset.). 
If we  we apply the definition of the Jones polynomial of virtual links to a classical links, 
it is the original Jones polynomial of the given classical link.   
The Jones polynomial of links in thickened surfaces and that of links in the 3-sphere have different properties  ((\cite{Kauffman1, Kauffman,Kauffmani}): 
\begin{figure}
\includegraphics[width=40mm]{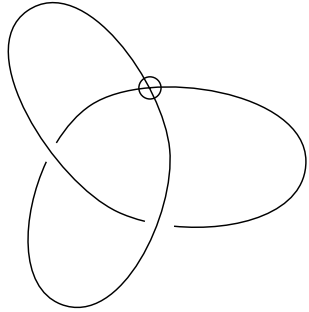}    
\caption{{\bf The Jones polynomial of this virtual knot is not 
that of any classical knot. }\label{vtre}}   
\end{figure}
There is a virtual 1-knot whose Jones polynomial is not that of any classical knot. 
An example is shown in Figure \ref{vtre}. 
A small circle placed around the crossing point as shown in Figure \ref{vtre} 
is called a virtual crossing. We will review it in \S\ref{psecv1k}.

A reason why virtual knot theory is important is as follows. 
It is an outstanding open question whether we can define the Jones polynomial 
in any 3-manifolds,  
although many other invariants are extended to 
the case of links in other 3-manifolds than the 3-sphere easily. 
Only one explicit partial answer 
is given  the Jones polynomial for links in thickened surfaces 
for now. It is given by Kauffman by using virtual links. 
The Jones polynomial of links in thickened surfaces and that of links in the 3-sphere have different properties as written above ((\cite{Kauffman1, Kauffman,Kauffmani}).  
Note that neither 
Reshetikhin and Turaev 
\cite[Theorem 3.3.3, page 560]{RT} nor 
Witten\cite{W} answers the above question. 
See Appendix. 

By using virtual links, we constructed a partial answer to the above question. 
Furthermore, 
Kauffman and Ogasa \cite{KOq} used virtual links, 
and introduced a new topological invariant 
$\upsilon_r$ 
of classical links as stated in the first paragraph.  
It is a main theme of this paper to discuss this invariant 
$\upsilon_r$.

In Part \ref{part2} of this paper we calculate the topological quantum invariants 
$\upsilon_r$ 
of  some examples explicitly. 
We conclude from our examples of explicit calculations 
that our invariant 
$\upsilon_r$ 
is new and strong enough to distinguish some classical knots from one another.  

Our main theorems are the following two claims.
\\

\def\Michi{
Our topological quantum knot invariants  
$\upsilon_r$ 
are strong enough to distinguish some classical knots from one another.}  

\h{\bf Theorem \ref{thmnontri}.}  
{\it \Michi}
\\

\def\Mni{
Our topological quantum invariants  
$\upsilon_r$ 
are strong enough to distinguish some 3-manifolds with boundary 
where the boundary is a disjoint union of two identical surfaces,  
from one another.}

\h{\bf Corollary \ref{cor3mfd}.}
{\it \Mni}

\vs

We review the new topological quantum invariants 
$\upsilon_r$ 
below.

When Jones  \cite{Jones} introduced the Jones polynomial,  
he \cite[page 360, \S10]{Jones} tried to define a 3-manifold invariant associated with the Jones polynomial, and succeeded in some cases. 
After that, 
Witten \cite{W} wrote a path integral for a 3-manifold invariant. 
 Reshetikhin and  Turaev \cite{RT} 
 defined a 3-manifold invariant 
 via surgery and quantum groups that one can view as  
 a mathematically rigorous definition  of the path integral. 
Kirby and  Melvin, and Lickorish and Kauffman and Lins
\cite{KM,Lickorish,Lickorishl,tl} continued this work.
Such 3-manifold invariants are called {\it quantum invariants $\tau_r$.}
These quantum  invariants  $\tau_r$ were defined for closed oriented 3-manifolds. 
In order to avoid confusion, we let 
$\upsilon_r$ 
denote our new topological quantum  invariant
and 
$\tau_r$  
the Reshetikhin-Turaev quantum invariant.

In \cite{KOq} Kauffman and Ogasa  
introduced topological quantum invariants 
$\upsilon_r$ 
of compact oriented  3-manifolds  with boundary 
where the boundary is a disjoint union of two identical surfaces. We review it in \S\ref{psecmth}.  
We explain how to use Kirby calculus for such manifolds, and
we use the diagrammatics of virtual knots and links to define these invariants.\\


Our invariants 
$\upsilon_r$ 
give new invariants of classical links in the 3-sphere in the following sense: Consider a link $L$ in $S^3$ of two components. \\
In this paper, the {\it complement} of a link means as follows: Take a tubular neighborhood $N(L)$ of $L$. $N(L)$ is the total space of the open $D^2$-bundle over $L$. The complement is $S^3-N(L)$. \\
The complement of a tubular neighborhood of $L$ 
is a manifold whose boundary consists in two copies of a torus. 
Our invariants apply to this case of bounded manifold 
and give new invariants of the given link of two components. 
We apply the same method and also obtain an invariant of 1-component links.  
See \S\ref{seccompL}.  
In this way, the theory of virtual links is used to construct 
new invariants of classical links in the 3-sphere.\\

It should be mentioned that the application of virtual knots to the calculation of these invariants is non-trivial and necessary. 
In \cite{DK},  Dye and Kauffman defined 
a quantum invariant for framed virtual links. 
In order to avoid confusion, we let 
$\varsigma_r$ 
denote the Dye-Kauffman quantum invariant. 
The Dye-Kauffman 
handling for Kirby calculus
and Temperley-Lieb Recoupling Theory for virtual link diagrams allows us 
to give specific formulas for our invariants for manifolds 
obtained by surgery on framed links embedded in a 
thickened surface. 
Just as the Jones polynomial can be calculated for links in thickened surfaces via virtual knot combinatorics, so can these surgery invariants be so calculated.
Note that in order to apply the virtual diagrammatic Kirby calculus, 
we need to set our diagrams so that the Roberts circumcision move 
$\mathcal O_3$ (Figure \ref{pfigO3})  is not needed. 
This we do by 
choosing a special surgical normalization as described below. One result of the normalization is that one cannot take any framed virtual diagram for our purposes, but any diagram can be modified so that the normalization is in effect.  
In this paper we review the definitions and frameworks of \cite{KOq}, and 
provide  specific calculations and applications.\\

In the sections to follow we address a number of issues. 
We show how to specify framings for the links in a thickened surface so that one can apply surgery.
We explain the results of Justin Roberts \cite{R} for surgery on three manifolds that are relevant and that apply for our use of Kirby Calculus. It should be noted that 
Robert's results use an extra move for his version of Kirby Calculus here denoted as $\mathcal O_3.$
 We show that the three manifolds that we construct can be chosen to 
have associated four manifolds that are simply connected, and that in this category, the topological types of these three manifolds are classified by just the first two of 
the moves $\mathcal O_1$ (Figure \ref{pfigO1}) and  $\mathcal O_2$ (Figure \ref{pfigO2}).
 Restricting ourselves to this category of three-manifolds, the first two moves correspond to the classical Kirby Calculus and to the generalized Kirby Calculus for virtual diagrams. \\

The $\mathcal O_1$ and  $\mathcal O_2$ moves 
do not change 
the Dye-Kauffman quantum invariants 
$\varsigma_r$ 
for framed virtual links.  
When Kauffman and Ogasa wrote the paper \cite{KOq}, 
it was open whether the $\mathcal O_3$ move changes the Dye-Kauffman quantum invariants 
$\varsigma_r$.  
Kauffman and Ogasa avoided the $\mathcal O_3$ move in order to introduce a topological invariant.   
Kauffman and Ogasa \cite{KOq} introduced a condition, {\it the  simple connectivity condition}, 
succeeded to avoid  the $\mathcal O_3$ move difficulty, and 
defined the new quantum invariants 
$\upsilon_r$ 
(see \S\ref{psec3mfd}). 

In this paper 
we proved the following about the above question. 
\\

\def\DK{
Let $F$ be a connected closed oriented surface. 
Let $M$ be a 3-manifold with boundary $F\amalg F$ 
with the boundary condition $\mathcal B$. 
Let $L^{fr}$ and $A^{fr}$ be framed links in $F\x[-1,1]$ with the symplectic basis condition $\mathcal F$ 
which represent $M$. 
Assume that $L^{fr}$ satisfies the simple connectivity condition but that 
 $A^{fr}$ does not satisfy it.  
 Then, in general,  
 each of the Dye-Kauffman quantum invariants 
 $\varsigma_r$ 
 of $L^{fr}$ 
 is not equal to that of  $A^{fr}$.}

\def\noO3{
In general, 
the $\mathcal O_3$ move on framed virtual links 
changes the Dye-Kauffman quantum invariants  
$\varsigma_r$  
for framed virtual links.}

\h{\bf Remark \ref{remnoO3}.}  
{\it \noO3}
\\

Remark \ref{remnoO3} implies the following. 
\\

\def\zettai{
In Definition $\ref{pdefdaiji}$,  the simple connectivity condition is necessary to define the 
%
our quantum invariants 
$\upsilon_r$.}

\h{\bf Remark \ref{remzettai}.}  
{\it \zettai}
\\

Note that we work only with closed oriented 3-manifolds that bound specific simply connected compact 4-manifolds, usually with these 4-manifolds corresponding to surgery instructions on a given link. Thus we concentrate on framed links that represent given 3-manifolds with boundary and that represent simply connected 4-manifolds.\\

In this way, we are able to apply Robert's results and make the connection between the topological types in a category of three manifolds and the Kirby Calculus classes of virtual link diagrams. With these connections in place, the paper ends with a description of the construction of 
the Witten-Reshetikhin-Turaev invariants $\tau_r$  that apply, via virtual Kirby Calculus, to our
category of three-manifolds. We obtain our new invariant 
$\upsilon_r$\\



We examine links that are embedded in thickened surfaces 
and obtain new invariants 
$\upsilon_r$
of three manifolds obtained by surgery 
on these thickened surfaces. 
One could take the viewpoint that the thickened surfaces are 
embedded in the three sphere and so also consider the three manifolds 
obtained by surgery on the links in the three sphere. 
These two points of view are distinct and give distinct invariants. 
This point of view for links in thickened surfaces 
is distinct 
from the usual point of view for the Reshetikhin Turaev invariants $\tau_r$, 
and our invariants 
$\upsilon_r$ 
give distinct results from these invariants.  
See \S\ref{secchiga}. 
\\




The invariants $\tau_r$ are defined only for framed classical links.
If  $L^{fr}$  is a framed classical link, 
the three invariants,  
$\tau_r$, 
$\varsigma_r$,  
and 
$\upsilon_r$,  
for 
$L^{fr}$ 
are equal.  
 The two invariants, 
 $\varsigma_r$ 
 and 
 $\upsilon_r$ 
 are defined for all framed virtual links. 
The invariants 
$\upsilon_r$ 
have different properties from the invariants 
$\varsigma_r$. 
It is a main theme of this paper. 
\\


As written above, the Dye-Kauffman invariants 
$\varsigma_r$ 
of framed virtual links 
do not produce a topological invariant of 3-manifolds 
if we do not impose the simple connectivity condition. 
In this paper we put an emphasis on this fact and we say that 
the Reshetikhin-Turaev invariants $\tau_r$ and our invariants 
$\upsilon_r$
are topological quantum invariants although we usually just say quantum invariants.

\part{The new quantum invariants 
$\upsilon_r$: 
Review of the definition}\label{part1}

\section{\bf 
Quantum invariants 
$\upsilon_r$ 
of 3-manifolds with boundary}\label{psecQI}

\begin{definition}\label{pdefB}
Let $M$ be a connected compact oriented 3-manifold with boundary. 
Let $\partial M=G\amalg H$, where $G$ and $-H$ are both orientation preserving diffeomorphic to 
a given closed oriented surface $F$ with genus $g$.  
Fix 
a symplectic basis 
$m^G_1,...,m^G_g$, 
(respectively, $m^H_1,...,m^H_g$),  
and longitudes, 
$l^G_1,...,l^G_g$, 
(respectively, $l^H_1,...,l^H_g$), 
for 
$G$ 
(respectively, $H$)
as usual.  
That is, the cohomology products of two basis elements are as follows:
$m^G_i\cdot l^G_i=+1$. 
$l^G_i\cdot m^G_i=-1$. The others are zero.  

Under these conditions, the 3-manifold $M$ is said to satisfy 
{\it boundary condition $\mathcal B$.}\\
\end{definition}

 We sometimes write $-F$ (respectively, $F\amalg -F$) as $F$ (respectively, $F\amalg F$) 
 when it is clear from the context.

We shall define 
topological quantum invariants of 
3-manifolds $M$ with boundary condition $\mathcal B$ below.\\

\noindent {\bf Remark:} For an oriented  manifold $X$, 
we sometimes write $-X$ as  $X$ when it is clear from the context.  \\

\section{
\bf Framed links in thickened surfaces}\label{psecFF}

It is well-known  that we can define the linking number for  2-component links 
in $\R^3$.  
On the other hand, in general, we cannot define it 
in the case of compact oriented 3-manifolds.  
However, we can define it in the case of thickened surfaces as below.



\begin{definition}\label{defSlk}
Let $(J,K)$ be a link in a thickened surface $F\x [-1,1]$, 
where $J$ and $K$ may be non-zero 1-cycles. 
We define the {\it linking number} lk$(J,K)$ of $J$ and $K$ 
as follows.
Let $Z$ stand for either $J$ or $K$.
The knot $Z$ together with a collection of circles in $F\x\{-1\}$ bounds 
a compact oriented surface $M_Z$ in $F\x[-1,1]$. 
%
%
Assume that 
$J$ 
(respectively, $K$)   
intersects $M_K$  
(respectively, $M_J$)    
transversely.   
Let $I(J, M_K)$ 
(respectively, $I(K, M_J)$) 
be the algebraic intersection number of $J$ and $M_K$ 
(respectively, $K$ and $M_J$). Note that $I(J, M_K)\neq I(K, M_J)$ in general. 
Let lk$(J,K)$ be \\ $\displaystyle\frac{1}{2}\{I(J, M_K)+I(K, M_J)\}$. 
\end{definition}

This is well-defined. It is proved by Reidemeister moves.

\begin{remark}\label{remFni}
(1) The linking number lk$(J,K)$ of a link $(J,K)$ in a thickened surface may be a half integer. 
Figure \ref{fignori} draws an example: The linking number of this link is 
$\frac{1}{2}$.   
	The places where a circle (inner and outer circles in the depiction of the torus surface) is cut indicate how the curves go from the upper part of the torus to the lower part (with respect to the projection directions chosen for this drawing). This discrimination allows us to indicate which crossings of the curves are actual weavings and which (dotted to solid curves) are artifacts of the projection.

\begin{figure}
\begin{center}  \includegraphics[width=800mm]{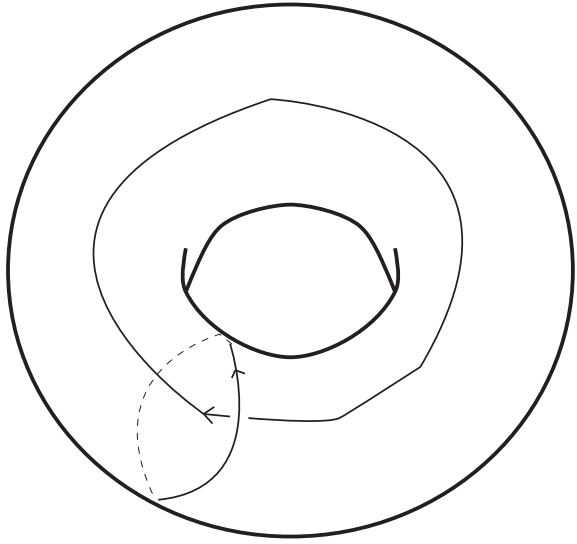}
\end{center}

\caption{{\bf 
A link in the thickened torus. The linking number is $\frac{1}{2}$.
}\label{fignori}}   
\end{figure}

If $F$ is the 2-sphere,  then $I(J, M_K)= I(K, M_J)$ holds and  lk$(J,K)$ is an integer. 
\\

\h(2) 
If $F$ is the 2-sphere, embed $F\x[-1,1]$ in the 3-sphere naturally. 
Regard a link $(J,K)$ in the thickened 2-sphere as a link in the 3-sphere. 
Then the well-known linking number of the link $(J,K)$ in the 3-sphere 
is equal to the linking number lk$(J,K)$ in Definition \ref{defSlk}. 


\end{remark}

An equivalent way to define linking numbers of links in thickened surfaces is as follows. 

\begin{definition}\label{defbe}
Let $(J,K)$ be a link in $F\x[-1,1]$. 
Make a projection of $(J,K)$ into $F\x\{-1\}$.    
Assume that the projection map is a self-transverse immersion. 
Just as we can take a diagram of a classical link in the plane, 
we can use such diagrams, which is the projection,  in the surface $F\x\{-1\}.$ 
We give each crossing point of $J$ and $K$ $+1$ or $-1$  
by using the orientation of $J$, that of $K$, and that of $F\x\{-1\}$.
The {\it linking number} lk$(J,K)$ is 
the half of the sum of the numbers at all crossing points.  
\end{definition}

This is well-defined. It is proved by Reidemeister moves.
The equivalence between two definitions above is also proved by  Reidemeister moves.










\bb
\h{\bf Framings.}  
Let $K$ be a knot in a compact oriented 3-manifold $M$. 
Recall that, in this paper, for a link $L$ in $M$, $N(L)$ is the total space of the open $D^2$-bundle over $L$. 
Let $\overline{N(K)}$ be the closure of the tubular neighborhood of $K$ in $M$. 
Take a knot $J$ in $\partial (\overline{N(K)})$ so that $J$ is homotopic to $K$ in $\overline{N(K)}$. 
We cannot define the linking number of $K$ and $J$ in general. 
If $M$ is a thickened surface, the linking number of $K$ and $J$ makes sense 
(Recall Definitions \ref{defSlk} and \ref{defbe}). 
Then it is an integer, not a half integer. 

\begin{remark}\label{remgold} 
For any 3-manifold $M$, 
we can always specify attaching maps of 4-dimensional 2-handles to $M$ by using a chart of $M$. 
%
%
However, we cannot determine the map by just choosing an integer for framing. 
The integer needs to be interpreted as a linking number of a curve on the boundary of a tubular neighborhood with the core of the solid torus. 
For this, it does suffice to have the surgeries on manifolds M of type $F\x[-1,1]$ where $F$ is a surface.
\end{remark}

Let $K$ be a knot in  a compact oriented 3-manifold $M$. 
Let $D^2\x B^2$ denote a 4-dimensional 2-handle such that 
$ D^2 \x \partial B^2$ is the attaching part. 
Let $O$ be the center of $D^2$, and $P$ a point in $\partial D^2$.  
We attach a 4-dimensional 2-handle $D^2\x B^2$ along a knot $K$ in $M$ 
so that $O\x \partial B^2$ coincides with $K$. 
Then $(O\x \partial B^2, P\x \partial B^2)$ is  a 2-component link in $M$.  
When we want to introduce framings associated with attaching 4-dimensional 2-handles, 
 we have to note the following fact.  
We cannot define the linking number of $O\x \partial B^2$ and  $P\x \partial B^2$ in $M$ in general.  
Therefore framings do not make sense in general. 
If $M$ is a thickened surface, framings make sense, and are always integers.

\begin{definition}\label{pdeffrlin} 
%
%
 A {\it framed link} $L^{fr}$ in a thickened surface is a link $L=(K_1,...,K_n)$ in a thickened surface such that each component $K_i$ is equipped with an integer in the sense described above so that this integer is a linking number. The integer is called the {\it framing} for $K_i$.
\end{definition}

Note that framings are always integers.

Assume that a component $K$ of framed link $L^{fr}$ in a thickened surface 
has a framing $n$. 
Take $\overline{N(K)}$ of $K$ in the thickened surface. 
 Embed a circle $C$ in $\partial (\overline{N(K)})$ which is homotopic to 
 $K$ in $\overline{N(K)}$ so that lk$(K,C)=n$. 
We attach a 4-dimensional 2-handle $D^2\x B^2$ along a knot $K$
so that $O\x \partial B^2$ coincides with $K$ 
and 
so that  $P\x \partial B^2$ coincides with $C$. 

Thus $L^{fr}$ 
represents, via framed surgery,  a compact oriented 3-manifold with boundary 
whose boundary is a disjoint union of the same two surfaces $F\amalg -F$.
It also represents a 4-manifold.

\bb
See Figure \ref{pfigtr} for an example.  
A framed link embedded in a thickened surface is drawn as 
a framed link in the surface which is the projection of  a thickened surface. 
Recall the explanation of drawing diagrams in Remark \ref{remFni}.(1). 

\begin{figure}
\begin{center}  \includegraphics[width=800mm]{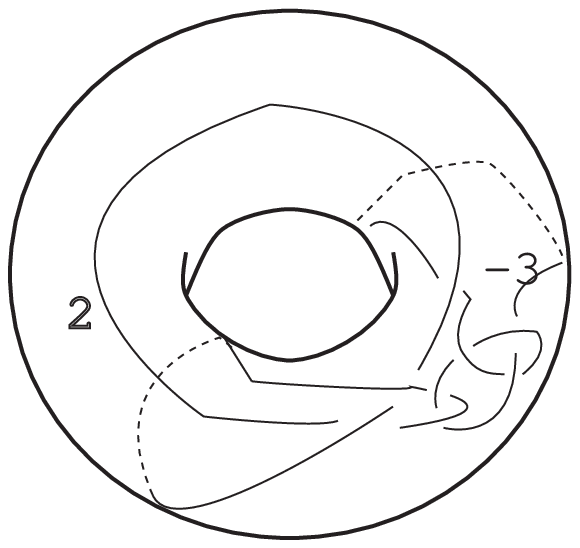}
\end{center}

\caption{{\bf 
A framed link in the thickened torus
}\label{pfigtr}}   
\end{figure}

\begin{definition}\label{pdefF}
Define $F\x[-1,1]$ with 
{\it the  symplectic basis condition $\mathcal F$} 
as follows. \\

Take $F\x[-1,1]$.
Fix a symplectic basis, 
$m^F_1,...,m^F_g$, 
and 
$l^F_1,...,l^F_g$, 
for 
$F$ as usual.
\\

Fix $m^F_i\x\{+1\}$ and  $l^F_i\x\{+1\}$ 
 (respectively,  $m^F_i\x\{-1\}$ and \\
 $l^F_i\x\{-1\}$)
in  $F\x\{+1\}$
(respectively,  $F\x\{-1\}$).  
\end{definition}

A framed link $L^{fr}$ in $F\x[-1,1]$ with 
the symplectic basis condition $\mathcal F$  
represents a compact oriented 3-manifold  
with the boundary $F\amalg F$  
with the boundary condition $\mathcal B$. 
It also represents a 4-manifold.\\
\\

\noindent{\bf Remark:} Note that another equivalent way to obtain framed links for the purpose of doing surgery on $F \times I$ is to use a generalized blackboard framing for a diagram drawn in the 
surface $F.$ 
Just as we can take a diagram of a classical link in the plane and regard it as a framed link by not using the first Reidemeister move and regarding the diagram itself as specifying a framing \cite{tl}, we can use such diagrams in the surface $F.$ 
In fact we can start with such a blackboard framed virtual link diagram (See \S\ref{psecv1k}), 
take the corresponding standard (abstract link diagram) construction producing a link diagram $L$ in a surface $F.$ The blackboard framing on the virtual diagram then induces a blackboard framing on the diagram in the surface. We will use this
association to show how the quantum link invariants we have previously defined for virtual link diagrams \cite{DK} become quantum invariants of actual three-manifolds via the constructions
in this paper.\\

Virtual links are represented by links in thickened surfaces.
We use these properties and define our quantum invariants.\\

\section{
\bf Framed link representations of 3-manifolds with non-vacuous boundary
}\label{psec3mfd}

Roberts \cite{R} generalized 
the result of Kirby\cite{Kirbyc} and  
the result of Fenn and Rourke\cite{FR}, and proved  
a theorem that is stronger than 
Theorem \ref{thm3mfd0} below.

\begin{thm}\label{thm3mfd0} {\bf  (This Theorem follows from Roberts' result \cite{R})}
Let $F$ be a closed oriented surface. 
Let $M$ be a compact oriented 3-manifold with boundary, 
whose boundary is $F\amalg F$,   
with the boundary condition $\mathcal B$.
Let $L_0$ and  $L_1$ be framed links in $F\x[-1,1]$ 
with the  symplectic basis condition $\mathcal F$,  
which represent $M$. 

Then 
$L_0$ and $L_1$  are related by the moves 
$\mathcal O_1$  in Figure $\ref{pfigO1}$,  
$\mathcal O_2$ in Figure $\ref{pfigO2}$, 
$\mathcal O_3$  in Figure $\ref{pfigO3}$, 
and framed isotopy in $F\x[-1,1]$.
\end{thm}

\h{\bf Remark:}
The $\mathcal O_1$ move  is carried out  in the 3-ball (Figure $\ref{pfigO1}$).   
The $\mathcal O_2$ move is carried out  in the genus two handle-body (Figure $\ref{pfigO2}$).
The $\mathcal O_3$   move is carried out  in the solid torus ,not the thickened torus. 
 (Figure $\ref{pfigO3}$).

\begin{figure}
\begin{center}\includegraphics[width=300mm]{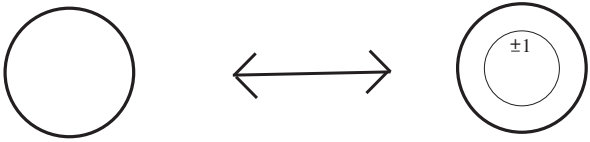}\end{center}
\caption{{\bf The operation $\mathcal O_1$ in the 3-ball}\label{pfigO1}}   
\end{figure}

\begin{figure}
\begin{center}\includegraphics[width=300mm]{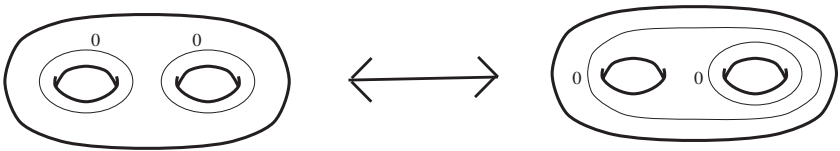}\end{center}
\caption{{\bf The operation $\mathcal O_2$ in the genus two handle-body
}\label{pfigO2}}   
\end{figure}

\begin{figure}
\begin{center}\includegraphics[width=300mm]{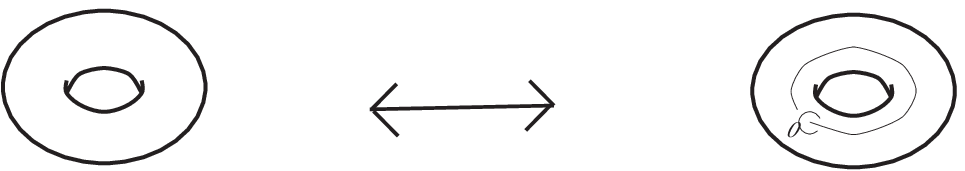}\end{center}
\caption{{\bf The operation $\mathcal O_3$ 
in the solid torus
}\label{pfigO3}}   
\end{figure}


\section{
\bf Framed links in thickened surfaces and 
simply connected 4-manifolds
}\label{psec3mfd}

We prove Theorem \ref{pthmbdr}
below, which is different from 
Theorem \ref{thm3mfd0}.   
This difference is very important. 
We use Theorem \ref{pthmbdr} and define our topological quantum invariant.

\begin{definition}\label{pdefS}
If a  framed link $L^{fr}$ in   $F\x[-1,1]$ represents a simply connected 4-manifold $V$,
we call $L^{fr}$ a framed link with 
{\it the simple-connectivity condition $\mathcal S$}. 
\end{definition}

By \cite{R}, 
 we have the following.

\begin{thm}\label{pthmbdr}
Let $M$ be a connected compact oriented 3-manifold  
with the boundary $F\amalg F$ 
with the boundary condition $\mathcal B$ $($in Definition $\ref{pdefB})$. 
Then $M$ is always described by a   framed link in   $F\x[-1,1]$ with 
the  symplectic basis condition $\mathcal F$ 
 $($in Definition $\ref{pdefF})$ 
and with 
the  simple-connectivity condition $\mathcal S$ 
$($in Definition $\ref{pdefS})$. 
 \end{thm}

\h{\bf Proof of Theorem \ref{pthmbdr}.}
Take a framed link $L^{fr}$ which represents $V$. 
Use the operation $\mathcal O_3$ in Figure \ref{pfigO3}, 
finitely many times,   
as shown in Figure \ref{pfigpi1}: 
Add a framed link to $L^{fr}$ as drawn in Figure \ref{pfigpi1}. 
Recall the explanation of drawing diagrams in Remark \ref{remFni}.(1).   

An example is drawn in Figure \ref{pfigexample}. 
\qed\\

We have also proved the following now.
\begin{thm}\label{thmtsuku}
Let $L^{fr}$ be a framed link 
$F\x[-1,1]$ with 
the  symplectic basis condition $\mathcal F$ $($in Definition $\ref{pdefF})$.  
 
Suppose that   $L^{fr}$ represents a connected compact oriented 3-manifold  
with the boundary $F\amalg F$ 
with the boundary condition $\mathcal B$ $($in Definition $\ref{pdefB})$. 

Then there is an explicit way to make  
the given framed link $L^{fr}$ 
into 
 a   framed link $A^{fr}$ in  $F\x[-1,1]$ with 
the  symplectic basis condition $\mathcal F$ 
and with 
the  simple-connectivity condition $\mathcal S$. 
Furthermore, 
 $L^{fr}$ 
is  a sub-framed link of $A^{fr}$. 
 \end{thm}

\begin{figure}
\centering
\includegraphics[width=160mm]{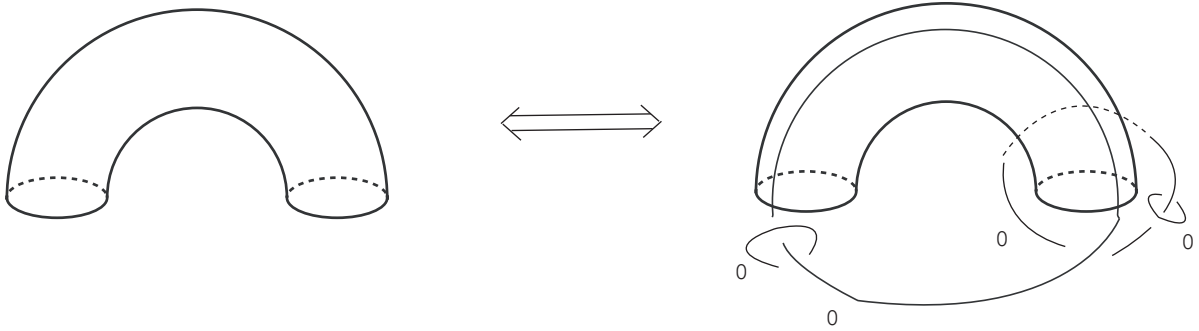}
\bb
\caption{{\bf 
Tow times of the $\mathcal O_3$ moves
}\label{pfigpi1}}   
\end{figure}

\begin{figure}
\centering
\includegraphics[width=160mm]{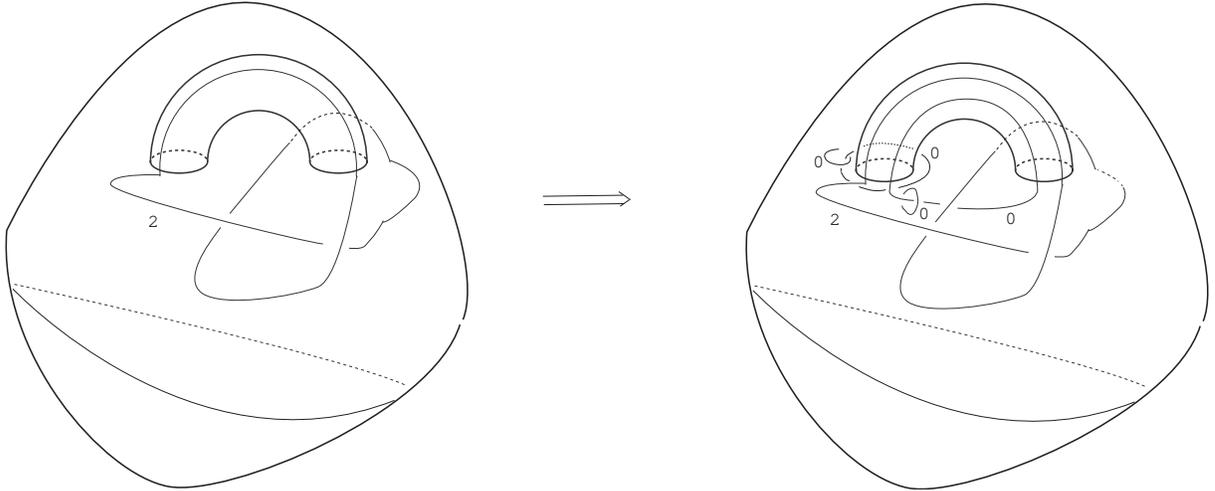}
\bb
\caption{{\bf 
Adding a framed link 
produces a new framed link 
with 
the  simple-connectivity condition $\mathcal S$
}\label{pfigexample}}   
\end{figure}

We generalize the results in 
Kirby, Fenn and  Rourke, and Roberts
\cite{FR, Kirbyc,R}, and we prove the following. \\

\begin{thm} \label{pthmKmove}
Let $F$ be a connected closed oriented surface.  
Let $M$ be a connected oriented compact 3-manifold  with boundary $F\amalg F$
with the boundary condition $\mathcal B$. 
Let $L^{fr}$ and ${L^{fr}}'$ be 
framed links in $F\x[-1,1]$ 
with 
the  symplectic basis condition $\mathcal F$ 
and 
with the  simple-connectivity condition $\mathcal S$,  
that  represent $M$. 
Then 
$L^{fr}$ is made from ${L^{fr}}'$ 
by a finite sequence of handle-slide, 
adding and removing the disjoint trivial knots with framing $\pm1$, 
that is, Kirby moves $($\cite{KM}$)$.
Note that under the hypothesis of this theorem we only use Roberts  moves 
$\mathcal O_1$ and $\mathcal O_2$.\\

Under these assumptions we have four manifolds $V$ and $V'$ as described 
in Definition \ref{pdefS}. 
Then \\
 $V\sharp^\alpha(S^2\x S^2)\sharp^\beta(S^2\widetilde\x S^2)
 \sharp^\gamma \C P^2\sharp^\delta \overline{\C P^2}$ 
is diffeomorphic to \\
 $V'\sharp^{\alpha'}(S^2\x S^2)\sharp^{\beta'}(S^2\widetilde\x S^2)
 \sharp^{\gamma'} \C P^2\sharp^{\delta'} \overline{\C P^2}$,
 where 
 $\alpha$, $\beta$, $\gamma$, $\delta$, 
 $\alpha'$, $\beta'$, $\gamma'$ and $\delta'$ 
  are non-negative integers.
 \end{thm}

\begin{remark}\label{pnotenonKmove}
Kirby moves mean the only $\mathcal O_1$ and $\mathcal O_2$ moves.

If we do not impose 
the  simple-connectivity condition $\mathcal S$ 
in Theorem \ref{pthmKmove}, 
$L^{fr}$ is not made from ${L^{fr}}'$ 
by a finite sequence of  Kirby moves in general. 
See Figure \ref{pfigbr}. 
Recall the explanation of drawing diagrams in Remark \ref{remFni}.(1).  
In the right figure of Figure \ref{pfigbr}, 
we draw a framed link in the torus which is the projection of  the thickened torus. 
The place where a circle is cut means 
which segments there goes over or down, as usual.
  The left figure of Figure \ref{pfigbr} represents  the empty framed link. 
The two framed links represent the same 3-manifold, 
but they are not Kirby move equivalent.
Note the difference between Figures \ref{pfigO3} and \ref{pfigbr}. 
In  Figure \ref{pfigO3}, we drew the solid torus.\\
\end{remark}

\begin{figure}
\centering
\includegraphics[width=500mm]{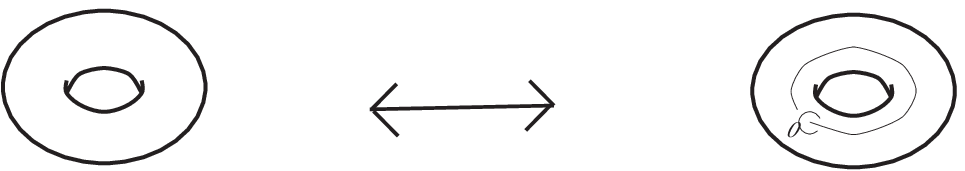}
\caption{{\bf 
Two framed links in the thickened torus 
}\label{pfigbr}}   
\end{figure}

\vs
The simple connectivity condition is important. 
Under the simple connectivity condition of framed links, 
we can carry out the $\mathcal O_3$ move 
by using the  $\mathcal O_1$ and  $\mathcal O_2$ moves.  
Figure \ref{fignoO3} is an example. 
\\

\begin{figure}
\centering
\includegraphics[width=600mm]{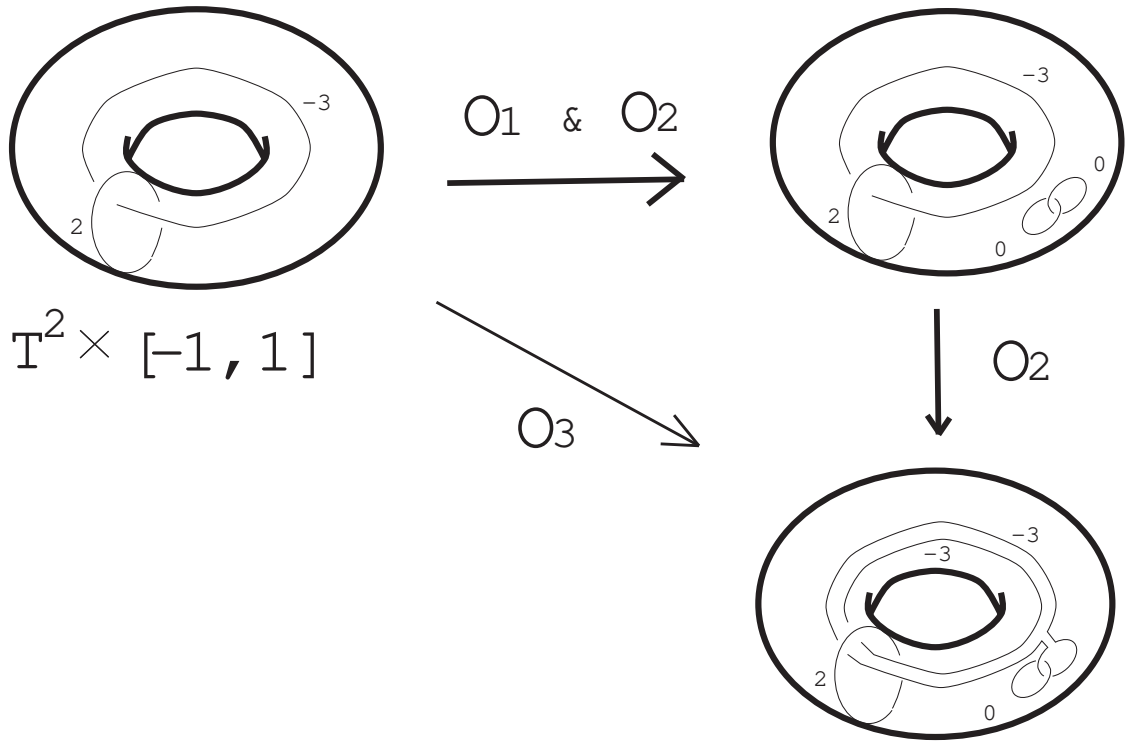}

\caption{{\bf 
An  $\mathcal O_3$ move realized by  $\mathcal O_1$ and  $\mathcal O_2$ moves
under the simple connectivity condition 
}\label{fignoO3}}   
\end{figure}

\section{
\bf Quantum invariants 
 $\varsigma_r$ 
 of framed virtual links: Outline 
}\label{psecvfr}

Kauffman \cite{Kauffman1,Kauffman,Kauffmani} describes and develops virtual links as a diagrammatic extension of classical links, and as a representation of links embedded in 
thickened surfaces.
The Jones polynomial of virtual links is defined in 
\cite{Kauffman1,Kauffman,Kauffmani}. 
See related open questions in \cite[\S4]{Org}. 

We can regard any framed link in $F\x[-1,1]$ 
as a framed virtual link.
See \cite{DK} for framed virtual links.
Note that the linking number of any pair, $K_i$ and $K_j$, is defined. The value is an integer or a half integer. Note that the framing is an integer. \\

Dye and Kauffman defined quantum invariants 
$\varsigma_r$  
of framed virtual links 
(See \S\ref{psecrevdef}).
If two framed links $L^{fr}$ and ${L^{fr}}'$ are changed into each other 
by a sequence of Kirby moves (\cite{KM}) and classical and virtual Reidemeister moves, 
each of the Dye-Kauffman quantum invariants 
$\varsigma_r$ 
of $L^{fr}$ is equivalent to that of ${L^{fr}}'$. \\

We use these invariants 
$\varsigma_r$ 
and,  introduce quantum invariants 
$\upsilon_r$ 
of 
3-manifolds with boundary in the following sections.\\

\begin{figure}
\includegraphics[width=60mm]{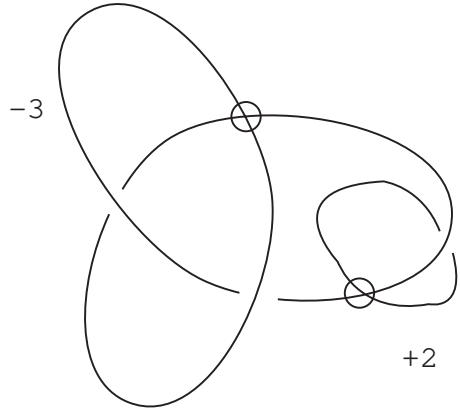}

\caption{{\bf 
A framed virtual link
}\label{pfigfrv}}   
\end{figure}

\section{ 
\bf Virtual knots and virtual links}\label{psecv1k}
\h
The theory of {\it virtual knot}s  
(\cite{Kauffman1,Kauffman,Kauffmani}) 
is 
 a generalization of classical knot theory, 
and  
studies the embeddings of circles in thickened 
oriented closed surfaces
modulo isotopies and orientation preserving diffeomorphisms
plus one-handle stabilization of the surfaces.   \\

By a one-handle stabilization, 
we mean a surgery on the surface that is performed on a curve 
in the complement of the link embedding and 
that either increases or decreases the genus of the surface. 
The reader should note that knots and links in thickened surfaces can be represented by diagrams on the surface in the same sense as link diagrams drawn in the plane or on the two-sphere. 
From this point of view, a one handle stabilization is obtained by cutting the surface along a curve in the complement of the link diagram and 
capping the two new boundary curves with disks, or taking two points on the surface in the link diagram complement and cutting out two disks, and then adding a tube between them. 
The main point about handle stabilization is that it allows the virtual knot to be eventually placed in a least genus surface in which it can be represented. 
A theorem of Kuperberg \cite{Kuperberg} asserts that such minimal representations are topologically unique. \\

\begin{figure}
\centering
\includegraphics[width=30mm]{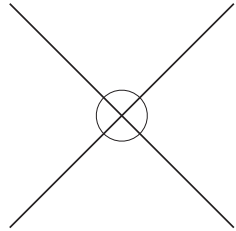}
\caption{{\bf Virtual crossing point}\label{pvcro0}}   
\end{figure}

Virtual knot theory has a diagrammatic formulation.
A {\it virtual knot} can be represented by a {\it virtual knot diagram} 
in $\R^2$ (respectively, $S^2$) 
containing a finite number of real crossings, and {\it virtual crossings} 
indicated by a small circle placed around the crossing point as shown in Figure \ref{pvcro0}.   
A virtual crossing is neither an over-crossing nor an under-crossing. A virtual crossing 
is a combinatorial structure keeping the information of the arcs of embedding
going around the handles of the thickened surface in the surface representation of the virtual link.\\

The moves on virtual knot diagrams
in $\R^2$ are generated by the usual Reidemeister moves plus the {\it detour move}.
The detour move allows a segment with a consecutive sequence of virtual crossings
to be excised and replaced by any other such a segment with a consecutive virtual
crossings, as shown in Figure \ref{pdetour}.    \\

Virtual 1-knot diagrams $\alpha$ and $\beta$ are changed into each other 
by a sequence of the usual Reidemeister moves  and detour moves 
if and only if 
$\alpha$ and $\beta$ are changed into each other 
by a sequence of all Reidemeister moves   
drawn in Figure \ref{pall-1}. 
\\

\begin{figure}
\includegraphics[width=110mm]{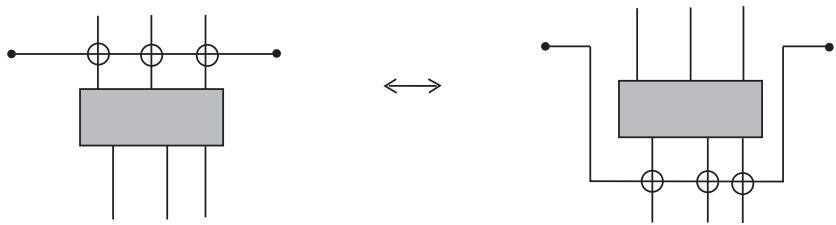}
\caption{{\bf An example of detour moves}\label{pdetour}}   
\end{figure}

\begin{figure}
\includegraphics[width=140mm]{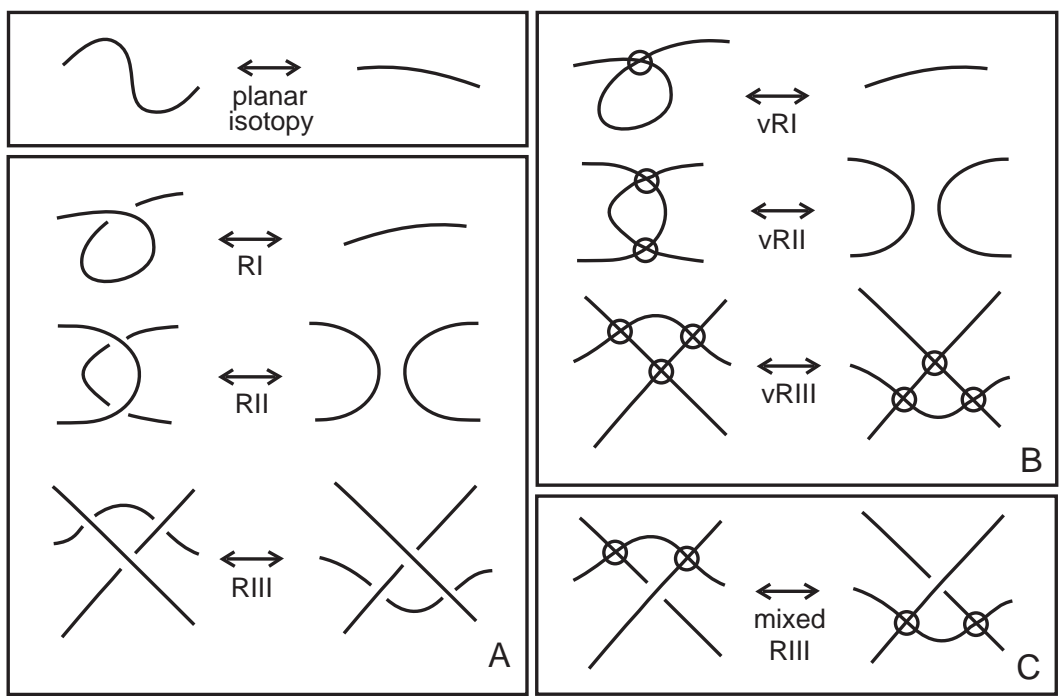}
\caption{{\bf All Reidemeister moves}   
\label{pall-1}}   
\end{figure}

Virtual knot and link diagrams that can be related to each other by a finite sequence of the
Reidemeister and detour moves are said to be
{\it virtually equivalent} or {\it virtually isotopic}. \\
The virtual isotopy class of a virtual knot diagram is called a {\it virtual knot}.\\

There is a one-to-one correspondence between the topological and the diagrammatic approach
to virtual knot theory. The following theorem providing the transition between the
two approaches is proved by abstract knot diagrams, see  \cite{Kauffman1,Kauffman, Kauffmani}.\\

\begin{thm}\label{pkihon} {\bf (\cite{Kauffman1,Kauffman, Kauffmani})} 
Two virtual link diagrams are virtually isotopic if and only if their surface embeddings are equivalent up to isotopy in the thickened surfaces, orientation preserving diffeomorphisms of the surfaces, and the  addition/removal of empty handles. 
\end{thm}

\begin{figure}
     \includegraphics[width=120mm]{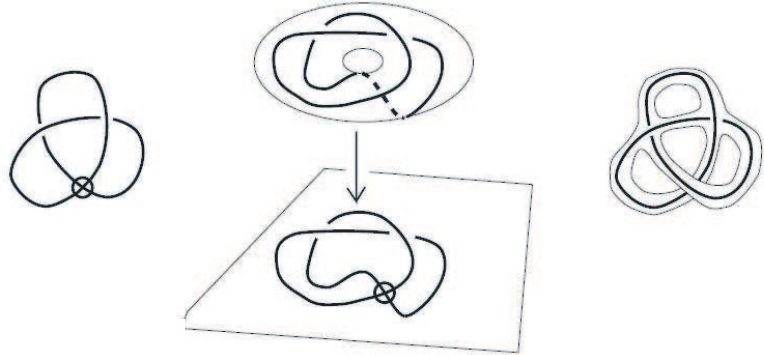}
\caption{{\bf 
How to make a representing surface from 
the tubular neighborhood of  a virtual knot diagram in $\R^2$
}\label{pvtube}}   
\bigbreak  
\end{figure}

\noindent{\bf Remark:}  
A handle is said to be {\it empty} if the knot diagram does not thread through the handle.  One way to say this more precisely is to model the addition of and removal of handles via the location of surgery curves in the surface that do not intersect the knot diagram.
Here, an oriented surface with a  link diagram using only classical crossings 
appears.
This surface is called a {\it representing surface}.   
In many figures of this paper, we use representing surfaces. 

In Figure \ref{pvtube}
we show an example of a way to make a representing surface from a virtual knot diagram. 
Take the tubular neighborhood of a virtual knot diagram in $\R^2$. 
Near a virtual crossing point, double the tubular neighborhood. 
Near a classical crossing point, keep    the tubular neighborhood and the classical crossing point. Thus we obtain a compact representing surface with non-vacuous boundary.  
 We may start with  a representing surface that is oriented and not closed, 
 and then  embed the surface in a closed oriented surface to obtain a new representing surface. 
 Taking representations of virtual knots up to such cutting (removal of exterior of neighborhood of the diagram in a given surface) and re-embedding, plus isotopy in the given surfaces, corresponds to a unique diagrammatic virtual knot type.
\\

\h{\bf The linking number.}
   
Recall the linking number of links in thickened surfaces in Definitions \ref{defSlk} and \ref{defbe}.  
If a link $(J,K)$ in $F\x[-1,1]$ is virtually equivalent to 
a link $(J',K')$ in $F'\x[-1,1]$, 
then we have ${\rm lk}(J,K)={\rm lk}(J',K')$.   
Therefore the following definition makes sense.  

\begin{definition}\label{defVlk0}
Let $(\mathcal J, \mathcal K)$ be 
a virtual link represented by 
a link $(J,K)$ in $F\x[-1,1]$. 
Define the {\it linking number} lk$(\mathcal J, \mathcal K)$ to be lk$(J,K)$. 
\end{definition} 

We have an alternative definition as below.

\begin{definition}\label{defVlk}
%
Let $(P,Q)$ be a virtual link diagram which represents a virtual link 
$(\mathcal P, \mathcal Q)$. 
Each classical crossing point of $P$ and $Q$ is oriented by 
the orientation of $P$, that of $Q$ and that of the plane. 
Assign to a positive classical crossing 
(respectively, a negative classical crossing, a virtual crossing) 
$+1$
(respectively, $-1$, $0$). 
Define lk$(P,Q)$ to be 
the half of the sum of  the number associated with each crossing point.  
Define the {\it linking number} lk$(\mathcal P, \mathcal Q)$ to be lk$(P,Q)$. 
\end{definition}

Examples are drawn in Figure \ref{figVHopf}. 
The linking number of the left virtual link is 
$-\frac{1}{2}$. 
The linking number of the right virtual link is 
$\frac{1}{2}$. 

\begin{figure}
\centering
\includegraphics[width=130mm]{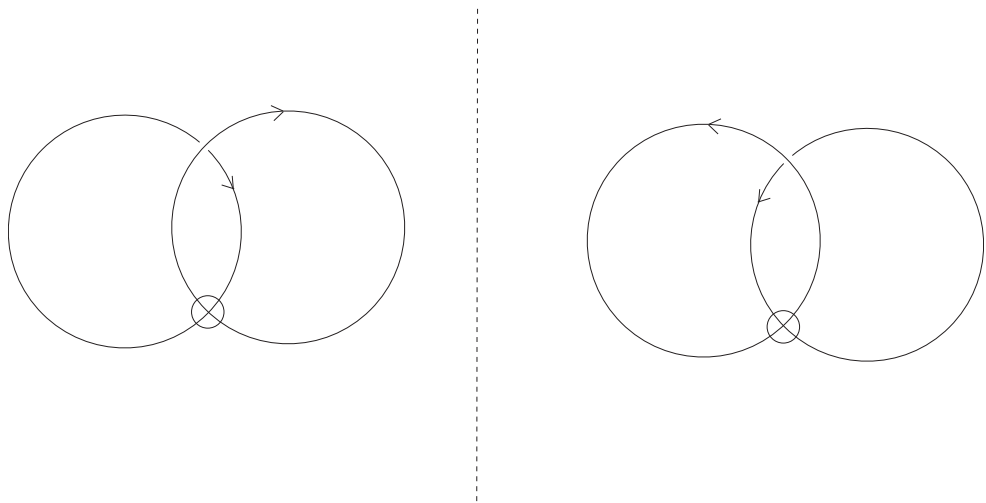}
\caption{{\bf Virtual Hopf link}\label{figVHopf}}   
\end{figure}


\vs
Manturov's textbook \cite{Mvbunrui} and 
its English translation \cite{IMvbunrui} by Ilyutko and Manturov 
are introduction to virtual knot theory.\\

\section{
\bf  Quantum invariants  
$\varsigma_r$ 
of framed virtual links: Review of Definition
}\label{psecrevdef}

\h
Dye and Kauffman \cite{DK}
 extended the definition of the Witten-Reshetikhin-Turaev invariant $\tau_r$
\cite{RT, RT2, W}  
to virtual link diagrams, and defined 
the Dye-Kauffman quantum invariants 
$\varsigma_r$
of framed virtual links. In this section $Z_K(r)$ denotes 
$\varsigma_r$ 
and we review the definition. See \cite{DK} for detail.\\

  First, we recall the 
definition of the Jones-Wenzl projector    
\\(q-symmetrizer) 
\cite{tl}.
We then define the colored Jones polynomial of a virtual link diagram. 
It is clear from this definition that two equivalent virtual knot diagrams have the same colored Jones polynomial. 
We  use these definitions to extend the Witten-Reshetikhin-Turaev invariant $\tau_r$
to virtual link diagrams. 
From this construction, we conclude that two virtual link diagrams, related by a sequence of framed Reidemeister moves and virtual Reidemeister moves, have 
the same value of the generalization 
$\varsigma_r$ 
of the Witten-Reshetikhin-Turaev invariant $\tau_r$. 
Finally, we prove that the 
generalization 
$\varsigma_r$ 
of the Witten-Reshetikhin-Turaev invariant $\tau_r$   
is unchanged by the virtual Kirby calculus. 
\\

To  form the \emph{n-cabling} of a virtual knot diagram, take $ n $ parallel copies of the virtual knot diagram. A single classical crossing becomes a pattern of 
$ n^2 $ classical crossings and a single virtual crossing becomes $ n^2 $ virtual crossings.\\

Let $ r $ be a fixed integer such that $ r \geq 2 $ and let 
\begin{equation*}
 A = e^{\frac{\pi i}{ 2r } }.
\end{equation*} 
Here is  a formula used in the construction of the Jones-Wenzl projector.  
\begin{equation*} \Delta_n = (-1)^n \frac{A^{2n + 2} - A^{-(2n+2)}}{A^2 - A^{-2}}. 
\end{equation*}
 Note that $ \Delta_1 = - (A^2 + A^{-2})$, the value assigned to a simple closed curve by the bracket polynomial.
 There will be an an analogous interpretation of $ \Delta_n $ which we will discuss later in this section.  \\
 
 We recall the definition of an n-tangle. Any two n-tangles can be multiplied by attaching the 
 bottom n strands of one n-tangle to the upper n strands of another n-tangle. We define an n-tangle to be \emph{elementary} if it contains no classical or virtual crossings. Note that the product of any two elementary tangles is elementary.
 Let $ I $ denote the identity n-tangle and let $ U_i $ such that $ i \in \lbrace 1,2, \dots n-1 
\rbrace $ denote the n-tangles shown in Figure \ref{pfig:tlgen}. \\


\begin{figure}
\begin{center}\includegraphics[width=60mm]{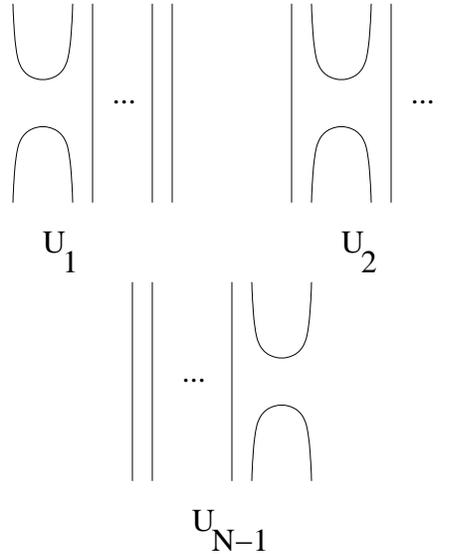}\end{center}
\caption{{\bf Tangles $U_i$ 
}\label{pfig:tlgen}}   
\end{figure}

By multiplying a finite set of $U_*$ as in  $ U_{i_1} U_{i_2} \dots U_{i_n} $, we can obtain any elementary n-tangle.
Formal sums of the elementary tangles over $ \mathbb{Z} [A, A^{-1} ] $ generate the $ n^{th} $ Temperly-Lieb algebra \cite{tl}.\\
 
We recall that the \emph{$ n^{th} $ Jones-Wenzl projector} is a certain sum of all elementary $ n $-tangles with coefficients in $ \mathbb{C}$ 
\cite{knotphys, tl}.
We denote the $ n^{th} $ Jones-Wenzl projector as 
$ T_n $. We indicate the presence of the Jones-Wenzl projector and the n-cabling by labeling the component of the knot diagram with $n$.\\

\begin{remark} 
There are different methods of indicating the presence of a Jones-Wenzl projector. 
In a virtual knot diagram, the presence of the $ n^{th} $ Jones-Wenzl projector is indicated by a box  with n  strands entering and n strands leaving the box. 
For n-cabled components of a virtual link diagram with an attached Jones-Wenzl projector, 
we indicate the cabling by labeling the component with $ n $ and the presence of the projector with a box. This notation can be simplified to the convention indicated in the definition of the the colored Jones polynomial. 
The choice of notation is dependent on the context.
\end{remark}

We construct the Jones-Wenzl projector recursively.  
The $1^{st} $ Jones-Wenzl projector consists of a single strand with coefficient $1$. 
There is exactly one 1-tangle with no classical or virtual crossings. 
The $q^{th} $ Jones-Wenzl projector is  constructed from the $ (q-1)^{th} $ and $ (q-2)^{th} $ Jones-Wenzl projectors as illustrated in figure \ref{pfig:asym}. \\


\begin{figure}
\centering
\includegraphics[width=120mm]{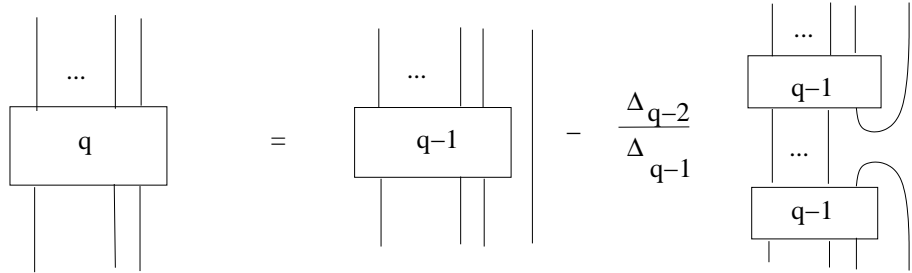}
\caption{{\bf 
$q^{th} $ Jones-Wenzl Projector
}\label{pfig:asym}}   
\end{figure}

We use this recursion to construct the $ 2^{nd} $ Jones-Wenzl projector as shown in figure \ref{pfig:2sym}.\\


\begin{figure}
\centering
\includegraphics[width=80mm]{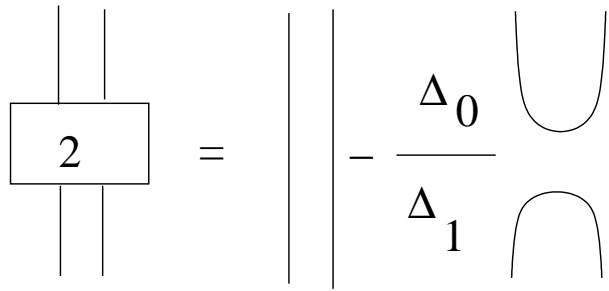}
\caption{{\bf 
$2^{nd} $ Jones-Wenzl Projector
}\label{pfig:2sym}}   
\end{figure}

We will refer to the Jones-Wenzl projector as the J-W projector for the remainder of this paper.\\

We review the properties of the J-W projector.
Recall that if $ T_n $ denotes the $ n^{th} $ J-W projector then


\begin{itemize}
\item[i)  ]  $T_n T_m = T_n$  \text{ for $ n \geq m $ }, \\
\item[ii) ] $T_n U_i = 0$  \text{ for all $ i $ }, \\
\item[iii)] \text{The bracket evaluation of the closure of } $T_n = \Delta_n$. 
\end{itemize}

\begin{remark}The combinatorial definition of the J-W projector  is given in \cite{tl}, p. 15.
Note that \cite{tl} provides a full discussion of all formulas given above. 
\end{remark} 

Let $ K $ be a virtual link diagram with components $ K_1, K_2 \dots K_n $. Fix an integer $ r \geq 2$ and let $ a_1, a_2 \dots a_n  \in  \lbrace 0, 1,2, \ldots r-2 \rbrace $. Let $ \bar{a} $ represent the vector $ ( a_1, a_2, \dots a_n ) $. Fix $ A= e^{\frac{\pi i}{2r}} $ and  $ d= -A^2 - A^{-2} $.
We denote the \emph{generalized $ \bar{a} $ colored Jones polynomial} of $ K $ as $ \langle K^{\bar{a}} \rangle $. 
To compute $ \langle K^{ \bar{a} } \rangle $, we
 cable the component $ K_i $ with $ a_i $ strands  and attach the $ a_i^{th} $ J-W projector  to  cabled component $ K_i $. We apply the Jones polynomial to the cabled diagram with attached J-W projectors. \\

The colored Jones polynomial is invariant under the framed Reidemeister moves and the virtual Reidemeister moves. This result is immediate, since the Jones polynomial is invariant under the framed Reidemeister moves and the virtual Reidemeister moves.\\

\begin{remark}The $ a $-colored Jones polynomial of the unknot is $ \Delta_{a} $. In other words, the 
 Jones polynomial of the closure of the $ a^{th} $ J-W projector is $ \Delta_a $.
\end{remark}

The \emph{generalized Witten-Reshetikhin-Turaev invariant} of a virtual link diagram is a sum of 
colored Jones polynomials. 
Let $ K $ be a virtual knot diagram with $ n $ components. Fix an integer $ r \geq 2 $. We denote the unnormalized Witten-Reshetikhin Turaev invariant of $ K $ as $ \langle  K^{ \omega}  \rangle $, which is shorthand for the following equation. 
\begin{equation}
\langle K^{\omega} \rangle = \underset{\bar{a} \in \lbrace 0,1,2, \dots r-2 \rbrace^n }{\sum}
\Delta_{a_1} \Delta_{a_2} \dots \Delta_{a_n} \langle K^{ \bar{a}} \rangle
\end{equation}

\begin{remark}
For the remainder of this paper, the Witten-Reshetikhin-Turaev invariant $\tau_r$  
will be referred to as the WRT.
\end{remark}

We define the matrix $ N $ in order to construct the normalized WRT \cite{tl}. 
Let $ N $ be the matrix defined as follows:

\begin{itemize}
\item[i) ] $N_{ij} =  lk (K_i,K_j) \text{ for } i \neq j $, \\

\item[ii) ] $N_{ii} = w(K_i)$. \\
\text{Then let} \\
        $b_{+} (K) = \text{ the number of positive eigenvalues of } N$, \\
        $b_{-} (K) = \text{ the number of negative eigenvalues of } N$,\\ 
    \text{and} \\
            $n(k) =  b_{+} (K) - b_{-} (K). $
\end{itemize}


The {\it normalized WRT} of a virtual link diagram $ K $ 
is denoted as $ Z_{K} (r) $.  
(This is the Dye-Kauffman quantum invariant 
$\varsigma_r$. 
See Remark \ref{remzet}. 
The Dye-Kauffman quantum invariants 
$\varsigma_r$ 
are defined 
only for framed virtual links. It is not defined for 3-manifolds. 
Kauffman and Ogasa \cite{KOq} used 
the Dye-Kauffman quantum invariants 
$\varsigma_r$, 
and   
introduced new topological quantum invariants 
$\upsilon_r$ 
for 3-manifolds.)

Let $ A = e^{ \frac{ \pi i }{ 2r}} $ and let $ |k|  $ denote the number of
components in the virtual link diagram $ K $. Then $ Z_K (r) $ is defined 
by the formula
\begin{equation*}
        Z_{K}(r) = \langle K^{ \omega}  \rangle \mu^{ | K | + 1}
        \alpha^{-n(K)}
\end{equation*}
where
\begin{align*}
  \mu &= \sqrt{ \frac{2}{r}} sin (\frac{\pi}{ r}) \\
\text{ and } \\
         \alpha &= (-i)^{r-2} e^{i \pi [ \frac{3(r-2)}{4r}]}.
\end{align*}
This normalization is chosen so that normalized WRT of the  
unknot with writhe zero is $ 1 $ and the normalization is invariant under the introduction and deletion of $ \pm 1 $ framed unknots.\\

Let $ \hat{U} $ be a $+1 $ framed unknot. We recall that $ \alpha = \mu \langle \hat{U}^{\omega} \rangle $ \cite{tl}, page 146. Since $ \hat{U} $ and $ K $ are disjoint in $ K \amalg \hat{U} $ then
 $ \langle (K \amalg \hat{U} )^{ \omega} \rangle = \langle K^{\omega} \rangle \langle \hat{U}^{\omega} \rangle $. We note that $ b_+ ( K \amalg \hat{U}) = b_+ (K) + 1 $,
 $ b_- ( K \amalg \hat{U}) = b_- (K) $, and $ | K \amalg \hat{U} | = |K|+1 $. We compute that
\begin{equation*}
Z_{K \amalg \hat{U} } (r) = \langle K^{ \omega} \rangle \langle \hat{U}^{\omega} \rangle \mu^{ |K|+2 } \alpha^{-n(K)-1}.
\end{equation*}
As a result, 
\begin{equation*} Z_{K \amalg \hat{U} } (r) = Z_K (r).
\end{equation*}\\


\begin{thm}\label{pinvar} 
Let $ K $ be a virtual link diagram then $Z_K(r)$ or 
$\varsigma_r$ 
is invariant under the framed Reidemeister moves, virtual Reidemeister moves, and the virtual Kirby calculus. 
\end{thm}

\h{\bf Remark:}  
 The virtual Kirby calculus means a sequence of only the $\mathcal O_1$ and $\mathcal O_2$ moves. 
Theorem \ref{pinvar}  never answers whether 
the $\mathcal O_3$ move changes the Dye-Kauffman invariants  
$\varsigma_r$ 
or not. 
Kauffman and Ogasa \cite{KOq} avoided answering this question, and
suceeded to introduce new topological invariants 
$\upsilon_r$
by using  the Dye-Kauffman invariants  
$\varsigma_r$. 
We review it in the following section \S\ref{psecmth}.

Furthermore, we prove in 
Proposition \ref{thmabsolutely} 
that 
the $\mathcal O_3$ move changes the Dye-Kauffman invariants  
$\varsigma_r$ 
in general.

\section{
\bf Our topological quantum invariants 
$\upsilon_r$ 
of 3-manifolds with boudary
}\label{psecmth}

\begin{definition}\label{pdefdaiji} 
Let $F$ be a connected closed oriented surface. 
Let $M$ be a connected oriented compact 3-manifold  
with the boundary $F\amalg F$  
with the boundary condition $\mathcal B$. 
Let $L^{fr}$ in $F\x[-1,1]$ with 
the  symplectic basis  condition $\mathcal F$ 
be  
a framed link with 
the  simple-connectivity condition $\mathcal S$ 
which represents $M$. 
Regard $L^{fr}$  as a framed virtual link. 
Define our {\it quantum invariants} 
$\upsilon_r$ 
of $M$ to be  
the Dye-Kauffman quantum invariants 
$\varsigma_r$ 
of 
the framed virtual link.
\end{definition}

By Theorem \ref{pthmKmove} and \S\ref{psecvfr}, we have the following.

\begin{mth} \label{pmthmain}
Definition $\ref{pdefdaiji}$ is well-defined, that is, 
each of the quantum invariants 
$\upsilon_r$ 
of $M$ is a topological invariant.
\end{mth}

%
%

\begin{remark}\label{remdiff}
Take a framed link in $F \times I$ with the symplectic basis condition. 
The value, 
$\upsilon_r(M)$, 
is invariant under diffeomorphisms of $F \times I$ by the property of virtual links.  
Therefore whatever symplectic basis for the symplectic basis condition $\mathcal F$ 
we take, we obtain the same value of our  invariants 
$\upsilon_r$. 
This is an important property for our invariants 
$\upsilon_r$. 
\end{remark}

\noindent {\bf Remark:} To actually apply our technique to a link $L$ in $F \times I$  we need to associate to the link $L$ embedded in $F \times I$ a framed virtual link diagram. The virtual Kirby class of this framed
virtual diagram will then be an invariant of the the three-manifold $M(L),$ and it is assumed that $L$ has been chosen so that the four manifold $W(L)$ is simply connected. 
As we have remarked, this condition of simple connectivity can be achieved by adding loops corresponding to the move $\mathcal O_3$ as illustrated in Figure~\ref{pfigpi1}.
If the link $L$ has originally been specified in $L \times I$ so that it satisfies the conditions of the Main Theorem, then one can obtain a virtual diagram for it, by taking a ribbon neighborhood of
a blackboard framed projection to $F$ and associating this with a virtual diagram in the standard way.\\

One can also start with a virtual diagram $K$ and associate an embedding  $L$ in $F \times I$ by the reverse process. However, the resulting $L$ may not satisfy the simple connectivity condition for the associated four-manifold. One way to insure this condition is to first associate a surface to $K$ by adding a handle to the plane at each virtual crossing. Then augment $K$ at each such handle to make sure that the simple connectivity condition is satisfied.
In Figure~\ref{pVA} we show how this augmentation of loops corresponds to an augmentation of a virtual diagram at a virtual crossing. Here we interpret the virtual crossing as corresponding
to a handle in the surface $F.$ We illustrate the augmentation at that handle an show how it corresponds to adding virtual curves to the given virtual diagram $K$ to form a virtual diagram $K'.$
If we start with a virtual diagram $K$ and apply this augmentation at each virtual crossing to form a virtual diagram $K',$ then resulting diagram $K'$ will represent a three manifold $M(K')$ that satisfies the simple connectivity condition. 
Thus the virtual Kirby class of this diagram $K'$  will be an invariant of the manifold $M(K').$ In this way we can create many examples for studying the results of this paper. In Figure~\ref{pVA1} we illustrate a specific example $K'$ whose invariants can be calculated. The reader interested in seeing the details of the calculation can consult 
\cite{tl,DK}, apply the above description of the invariants and work out the expansion of the invariants for the link $K'$ in Figure~\ref{pVA1} .

\begin{figure}
\centering
\includegraphics[width=80mm]{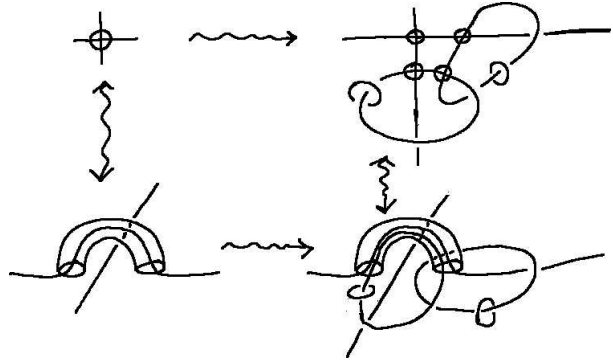}
\caption{{\bf Virtual Augmentation}
\label{pVA}}   
\end{figure}

\begin{figure}
\centering
\includegraphics[width=80mm]{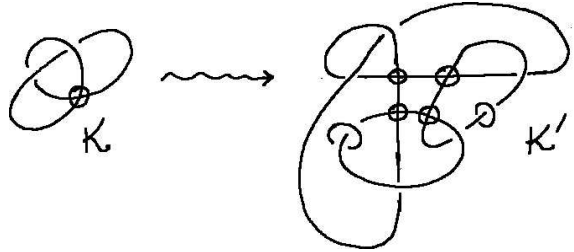}
\caption{{\bf Virtual Augmentation Example}
\label{pVA1}}   
\end{figure}

\vs

\section{
\bf Our topological quantum invariants 
$\upsilon_r$ 
of 
classical knots in the 3-sphere
}\label{secknot}

We define topological quantum invariants 
$\upsilon_r$ 
of knots in the 3-sphere. 

We make a knot $K$ in $S^3$ into a 2-component link $L=(K,J)$ in $S^3$ as follows: 
$J$ is the trivial knot. 
There is an embedded 2-disc in $S^3$ 
that $J$ bounds and that $K$ intersects geometrically once. 
Thus the linking number of $K$ and $J$ is one when we give an orientation to $L$. 
Then we say that $L=(K,J)$ is the {\it ring-hooked knot} of $K$.
Note that the ring-hooked knot of $K$ is determined by $K$ uniquely. 
When  $K$ has a well-known name, e.g., the trefoil knot, 
we abbreviate `the ring-hooked knot of the trefoil knot' with `the ring-hooked trefoil knot'.
The ring-hooked trivial knot is the Hopf link. 

The complement of a ring-hooked knot $(K,J)$ is a compact oriented 3-manifold 
whose boundary is the disjoint union of two tori. 
We can define our topological quantum invariants 
$\upsilon_r$  
for the complement if we induce the boundary condition $\mathcal B$.
We put a symplectic basis for $\partial(S^3-N(K\amalg J))$ 
in the next two paragraphs.
Call this basis the {\it standard basis}. 

We put a basis $(u_1, v_1)$ for $K\x \partial D_2$ as follows: 
$u_1$ is defined by the meridian of $K$. 
$v_1$ is defined by 
a circle $C$ embedded in $K\x \partial D_2$ 
such that $C$ is homotopic to $K$ in $K\x D^2$ and 
such that the linking number of $K$ and $C$ in the 3-sphere is zero.

We put a basis $(u_2, v_2)$ for $J\x \partial D_2$ as follows: 
$u_2$ is defined by 
a circle $E$ embedded in $J\x \partial D_2$ 
such that $E$ is homotopic to $J$ in $J\x D^2$ and 
such that the linking number of $J$ and $E$  in the 3-sphere is zero. 
$v_2$ is defined by the meridian of $J$.

 \h Remark: We have $[u_1]=[u_2]$ and $[v_1]=[v_2]$ in $H_1(S^3-N(K\amalg J);\Z)$, 
 where $u_1$, $u_2$, $v_1$, and $v_2$ represent circles.  

We define a topological quantum invariants 
$\upsilon_r$ 
of each knot $K$ 
to be our topological quantum invariants 
$\upsilon_r$ 
of 
the complement of the ring-hooked knot of $K$ 
with the above symplectic basis. 

By our construction, 
$\upsilon_r$ 
is a topological  invariant of knots in the 3-sphere. 

\bb

\part{\bf 
Calculation}\label{part2}

\section{\bf 
The simple connectivity condition is necessary to define our quantum invariants 
$\upsilon_r$ 
}\label{secabso}

\begin{figure}
\includegraphics[width=90mm]{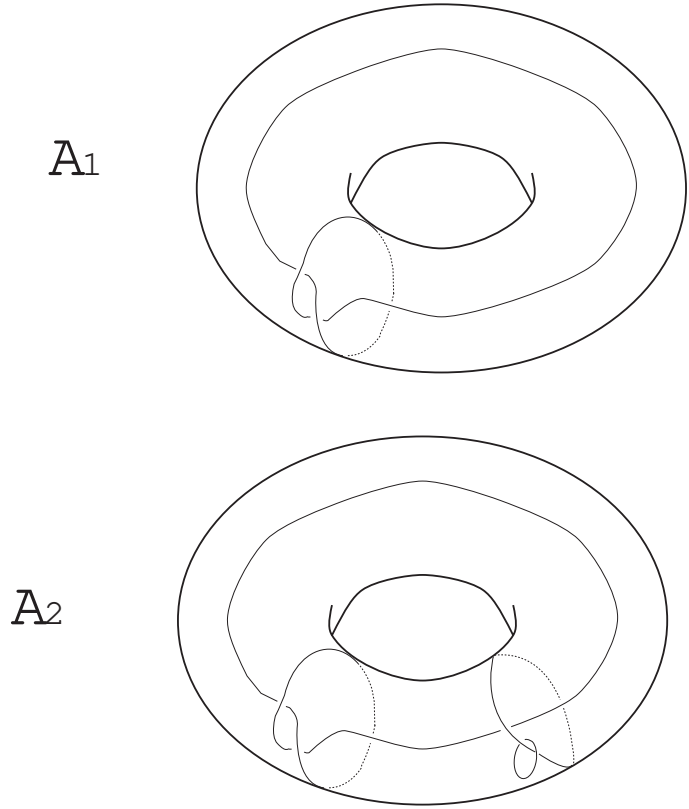}
\caption{{\bf 
Framed links in thickened tori. All framings are zero. 
}\label{figcheck}}   
\end{figure}


 Figure \ref{figcheck} 
shows framed links $A_1$ and $A_2$ in the thickened torus.
$A_2$ is obtained from $A_1$ by one $\mathcal O_3$ move. 
Therefore $A_1$ and $A_2$ represent the same 3-manifold with boundary,  
with the same boundary condition $\mathcal B$ $($Definition $\ref{pdefB})$.
 We have the following. 

\begin{proposition}\label{thmabsolutely}
In general, each of the Dye-Kauffman quantum invariants 
$\varsigma_r$ 
of $A_1$
 is not equal to that of $A_2$ 
although $A_1$ and $A_2$  represent the same 3-manifold with boundary 
with the same boundary condition $\mathcal B$. 
\end{proposition}

\h{\bf Proof of Proposition \ref{thmabsolutely}.}
Calculation. See Table \ref{tab:results}.  
\qed\\






Therefore we must concentrate on only framed links which satisfy the simple connectivity condition.
We have the following observations. 
\\

\begin{remark}\label{remzet}
{\it \DK}
\end{remark}

\h{\bf Remark:}  
Our quantum invariants 
$\upsilon_r$ 
are not defined for $A^{fr}$  
because  $A^{fr}$ does not satisfy the simple connectivity condition. 
The Dye-Kauffman invariants 
$\varsigma_r$ 
are defined for $A^{fr}$ and 
for all links in thickened surfaces because it is defined for all framed virtual links.  

\begin{remark}\label{remnoO3}
{\it \noO3}
\end{remark}

\begin{remark}\label{remzettai}
{\it \zettai}
\end{remark}

\vs

Let $L$ and $L'$ be framed links in a thickened surface. 
Assume that $L$ is obtained from $L'$ by one $\mathcal O_3$ move. 

If $L$ and $L'$ satisfy the simple connectivity condition, 
$L$ is obtained from $L'$ by a sequence of only the $\mathcal O_1$ and $\mathcal O_2$ moves. 
Recall Theorem \ref{pthmKmove}. 
Therefore each of the Dye-Kauffman quantum invariants 
$\varsigma_r$ 
of $L$ and that of $L'$ are equivalent, 

Remark \ref{remnoO3} claims as follows: 
Assume the condition $(\ast)$ that one of $L$ and $L'$ does not satisfy the simple connectivity condition. 
In general, 
each of the Dye-Kauffman quantum invariants 
$\varsigma_r$ 
of $L$ is not equivalent to that of $L'$.

However, under this condition $(\ast)$, 
it can happen that  
the $\mathcal O_3$ move does not 
change the Dye-Kauffman quantum invariants 
$\varsigma_r$.
See Figure \ref{figcheck00Y}. 
\begin{figure}
\includegraphics[width=130mm]{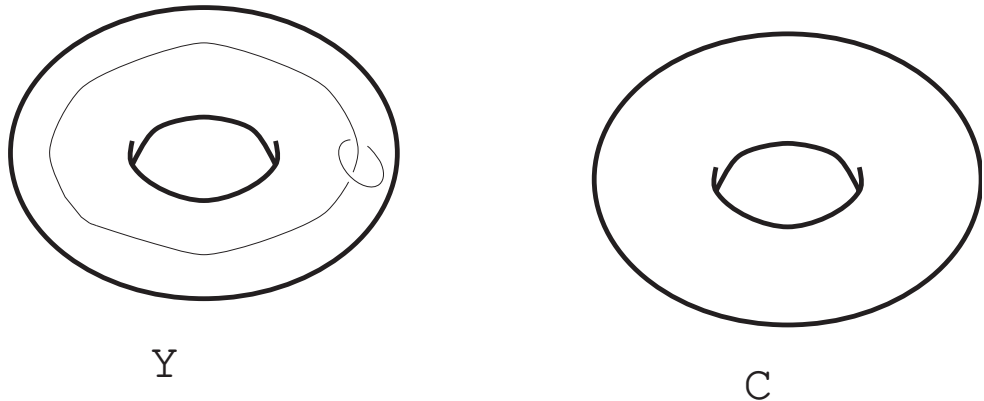}
\caption{{\bf  
Framed links in thickened tori. All framings are zero. 
}\label{figcheck00Y}}   
\end{figure}
$C$ is the empty framed link. 
$Y$ is a 4-component framed link such that all framings are zero in the thickened torus.  
 $Y$ and $C$ represent the same 3-manifold. 
 $Y$ and $C$ represent different 4-manifolds 
such that their fundamental groups are different. 
Therefore $Y$ is not obtained from $C$ by 
a sequence of the $\mathcal O_1$ and the $\mathcal O_2$ moves.  
However, each of the Dye-Kauffman quantum invariants 
$\varsigma_r$ 
of $C$ and 
that of $Y$ are the same. \\
{\it Reason}:   
The virtual framed link representation of $Y$ is the classical Hopf link such that both framings are zero.
Therefore it is changed into the empty framed link 
by the $\mathcal O_1$ and the $\mathcal O_2$ moves. 
\\

See also the question in the last part of \S\ref{seckekka}.

\section{
\bf Framed links in complements of knots}\label{seccompK}
As written in \S\ref{psecintro}, 
in this paper, the {\it complement} of a link means as follows: 
Take a tubular neighborhood $N(L)$ of $L$. $N(L)$ is the open $D^2$-bundle over $L$. 
The complement is defined to be $S^3-N(L)$.

See Figure \ref{fighokuK2}. 
The left upper figure is the complement of the trivial knot $K$, 
which includes a knot $J$. 
$J$ is null-homologous in $S^3-K$. 
Examine the top two diagrams in Figure \ref{fighokuK2}. 
We attach a 4-dimensional 2-handle to the complement of $K$ along $J$ 
as follows. 
Since $J$ is null-homologous in the complement of $K$, 
the framing of $J$ makes sense. 
Let the framing  be $+1$.

The result of this surgery 
is drawn in the right figure. 
It is the complement of the right-handed trefoil knot $K'$. 
\\

\h{\bf Remark:} The complement of $K$ is the solid torus. 
In general, we cannot define the linking number of 2-component links 
or the framing on knots 
in the solid torus. Recall \S\ref{psecFF}.   

However, the complement of $K$ is embedded in the 3-sphere as drawn in  Figure \ref{fighokuK2}.   
The framing of $J$ in the complement of $K$ is specified by using the 3-sphere. 
Recall Remark \ref{remgold}.  
In Figure \ref{fighokuK2}, 
these two ways of the definitions of framing coincide.
{\it Reason}: Since  $J$ is null-homologous in the complement of $K$, 
there is a compact oriented surface $F$ in the complement of $K$ whose boundary is $J$. 
$J$ is null-homologous in the 3-sphere and $F$ is embedded in the 3-sphere.
In both cases, the framing can be defined by using $F$. 
\bb


As for the rest of the diagrams in Figure  \ref{fighokuK2}, 
the lower three figures are isotopic to the left upper figure. 
Thus one could have begun with the lowest diagram and 
pointed out that a surgery on 
the framed knot with framing $+1$
would produce the complement of the  right-handed trefoil knot in the three sphere. We will use this kind of surgery in discussing links in thickened surfaces in the following	sections.
\\

%





\begin{figure}
\begin{center}\includegraphics[width=110mm]{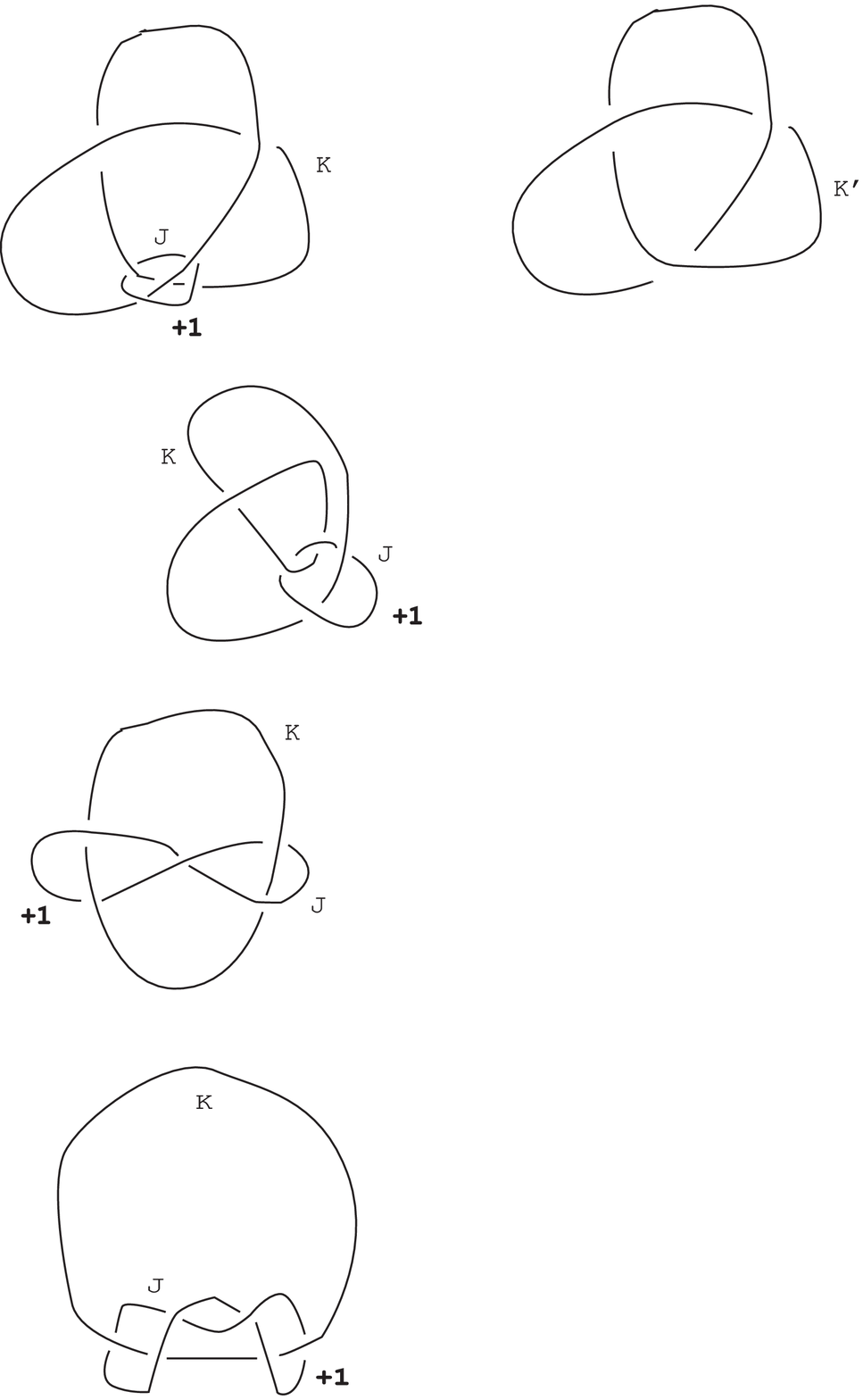}\end{center}
\caption{{\bf 
A `framed link' in the complement of $K$ 
and the result of the surgery
}\label{fighokuK2}}   
\end{figure}

\section{
\bf Framed links in complements of 2-component links
}\label{seccompL}


We calculate examples of our topological quantum invariants 
$\upsilon_r$  
of knots in the 3-sphere defined in \S\ref{secknot}. 



The first example is the right-handed trefoil knot. 
We make the ring-hooked  right-handed trefoil knot $(K'_1, K'_2)$ 
with the standard basis,   $(u_1, v_1)$ and $(u_2, v_2)$, 
(defined in \S\ref{secknot})
as drawn in the right upper figure of Figure \ref{fighokuL1}.  
We have that 
$\partial(S^3-N(K'_1\amalg K'_2))$ is a disjoint union of two tori.  
In order to calculate  our topological quantum invariants 
$\upsilon_r$, 
we  construct 
a framed link in 
the thickened torus with the symplectic basis condition $\mathcal F$
 which represents $S^3-N(K'_1\amalg K'_2)$ 
with the standard basis 
as below.

Take the ring-hooked trivial knot $(K_1,K_2)$. 
with the standard basis
 $(p_1, q_1)$ in $\partial N(K_1)$ 
 and 
$(p_2, q_2)$ in $\partial N(K_2)$. 
Note that 
the complement  of the ring-hooked trivial knot $(K_1,K_2)$ is the thickened torus. 
The complement with the standard basis is 
the thickened torus with the simplectic basis condistion $\mathcal F$. 
Recall Remark \ref{remdiff}: 
Each of our topological quantum invariants 
$\upsilon_r$ 
for any of the symplectic basis for the symplectic basis condition $\mathcal F$ is the same.

The left upper figure in Figure \ref{fighokuL1} 
is the complement of the Hopf link $(K_1, K_2)$ 
which includes the framed knot $K$ with framing $+1$.

Note that framings make sense in the case of links in thickened surfaces as we review in \S\ref{psecFF}. 
%
%
%
Since $K$ is null-homologous, we can also define a framing on $K$ 
by using the fact that the $K$ is null-homologous. 
The framing is also $+1$.

We carry out $+1$ surgery along the framed knot on the complement of $(K_1, K_2)$. 
The result is drawn in the right figure. 
It is the complement of 
the ring-hooked right-handed trefoil knot 
$(K'_1, K'_2)$.  
Since $K$ is null-homologous, 
the basis 
$(p_1, q_1)$ 
(respectvely, $(p_2,q_2)$) 
for 
$K_1$ 
(respectvely, $K_2$)   
is changed into 
the basis 
$(u_1, v_1)$ 
(respectvely, $(u_2,v_2)$) 
for 
$K'_1$ 
(respectvely, $K'_2$)  
after the surgery along $J$ with framing $+1$. 
{\it Reason}: There is a compact oriented surface 
$F_p$ 
(respectively, $F_q$)
with boundary 
in $S^3-N(K_1\amalg K_2)$ 
whose boundary is 
$p_1$ and $p_2$ 
(respectively, $q_1$ and $q_2$).  
Note that $F_p$ and $F_q$ intersect.
Since $J$ is null-homologous in 
in $S^3-N(K_1\amalg K_2)$, 
$J$ never intersects  
$F_p$ 
(respectively, $F_q$)
algebraically.

The left middle figure is obtained from the left upper figure by an isotopy. 

In the left lower figure, 
we draw a framed 5-component link in the complement of $(K_1, K_2)$: 
One component has framing $+1$ and the others have framings $0$. 
It is obtained from 
the framed link in the left middle figure 
by two times of the $\mathcal O_3$ move.   
Of course, both framed links represent the same 3-manifold.

Remark:  
The framed link in the left lower figure satisfies the simple connectivity condition. 
The one in the left middle figure does not.
\\






Remark: Figure \ref{fighokuL1} draws a link $(K_1,K_2)$ in the 3-sphere.
The complement of $(K_1,K_2)$ in the 3-sphere is the thickened torus. 
The framing of the knot 
by using the thickened torus, or the complement,  is $+1$ (Recall \S\ref{psecFF}). 

The complement of $(K_1,K_2)$ is embedded in the 3-sphere.
We also specify the framing by using the 3-sphere. It is also $+1$.

\begin{remark}\label{remreirei} 
Let $Q$ be the complement of the Hopf link $(K_1, K_2)$ as drawn in Figure \ref{figreirei}.
Here, note that 
the complement is the thickened torus   
and that 
$Q$ is a submanifold of $S^3$. 

Take a framed knot $J$ in $Q$ whose framing is defined by using the thickened torus. 
We can also defined a framing on $J$ by using $S^3$ 
because 
$Q$ is a submanifold of $S^3$.

In general, both methods do not give the same framing. 
 See Figure \ref{figreirei}. 
The boundary $\partial Q$ is two tori.   
Take a collection $X$ of circles in one of those two tori
such that $X$ and $J$ are homologous in $Q$.  
Since $Q$ is a submanifold of the 3-sphere, 
 the linking number $\rho$ of $(J,X)$ in the 3-sphere is determined uniquely by using the 3-sphere.  
Both methods give the same framing if and only if 
$\rho$ 
is $0$. 
In Figure \ref{figreirei}, we have $\rho\neq0.$  

Framings in  figures of
this paper 
are defined by using the thickened surfaces. 
\end{remark}
\bb

The framed link in the left lower figure in Figure \ref{fighokuL1} 
is isotopic to that in Figure \ref{figdrawkae}.   
These framed links, 
that in Figure \ref{figdraw},  
that in Figure \ref{figdrawY}, and 
that in  Figure \ref{figdrawX}  
represent 
the same 3-manifold with the same boudary condition $\mathcal B$
and 
the same 4-manifold 
because they are changed into each other by 
handle slides. 
See Figure \ref{figdrawX} and imagine a handle slide, 
which changes Figure \ref{figdraw} into Figure \ref{figdrawY} 
(respecctively,  Figure \ref{figdrawkae} to Figure \ref{figdraw}).










\begin{figure}
\begin{center}\includegraphics[width=130mm]{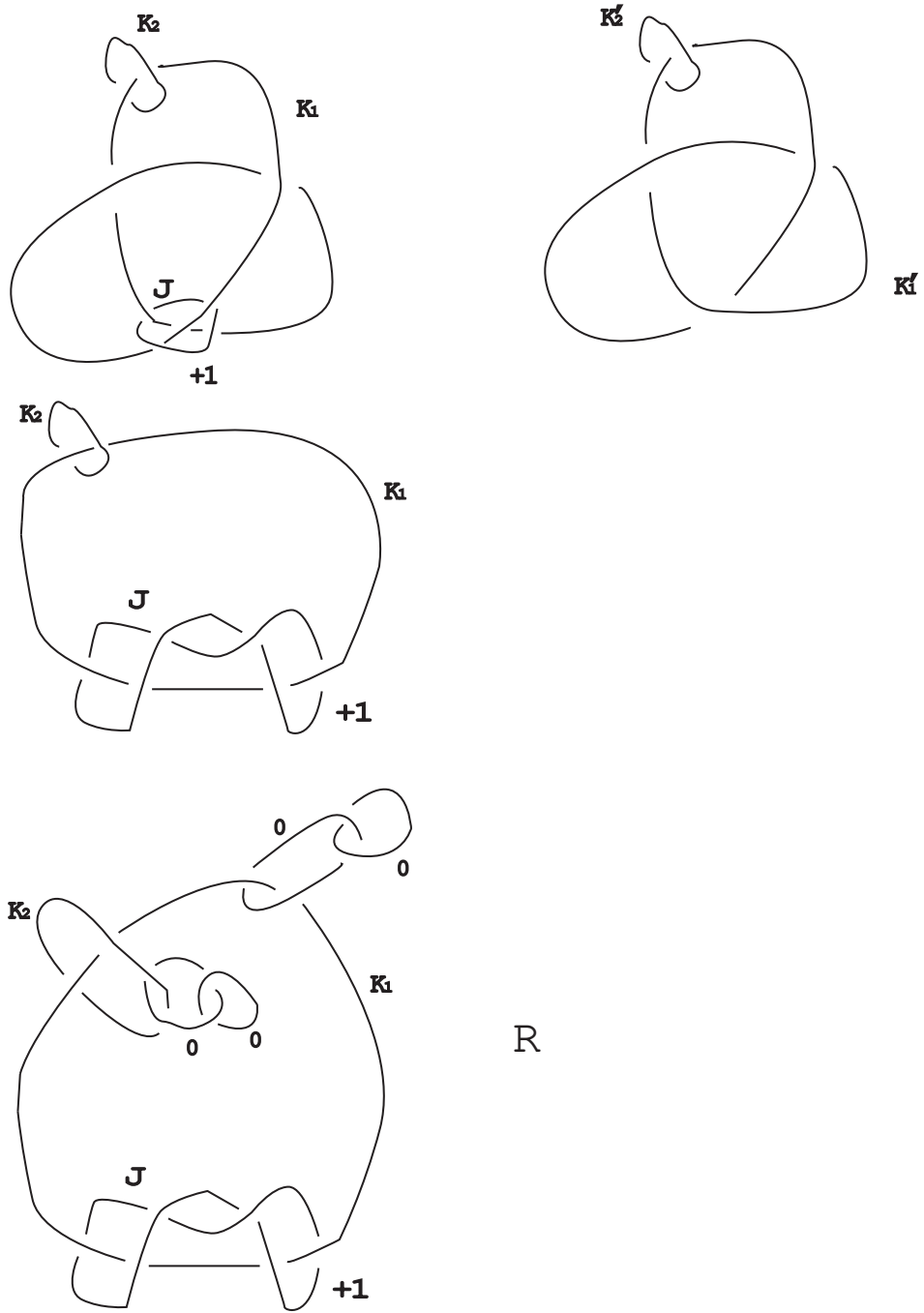}\end{center}

\caption{{\bf 
A `framed link' in the complement of $(K_1,K_2)$ 
and the result of the surgery
}\label{fighokuL1}}   
\end{figure}

\begin{figure}
\includegraphics[width=70mm]{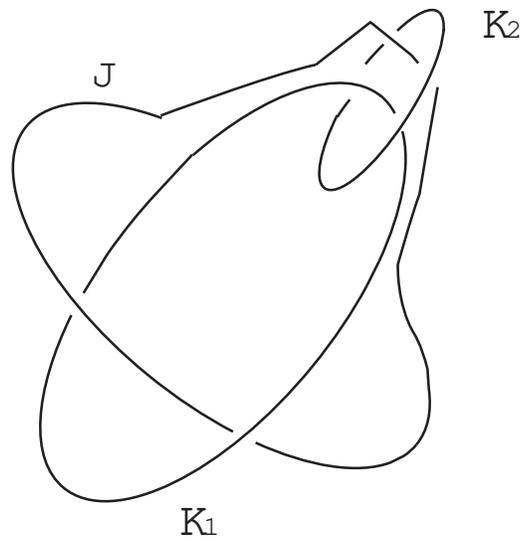}
\caption{{\bf 
An example for Remark \ref{remreirei}.
}\label{figreirei}}   
\end{figure}

\begin{figure}
\begin{center}\includegraphics[width=130mm]{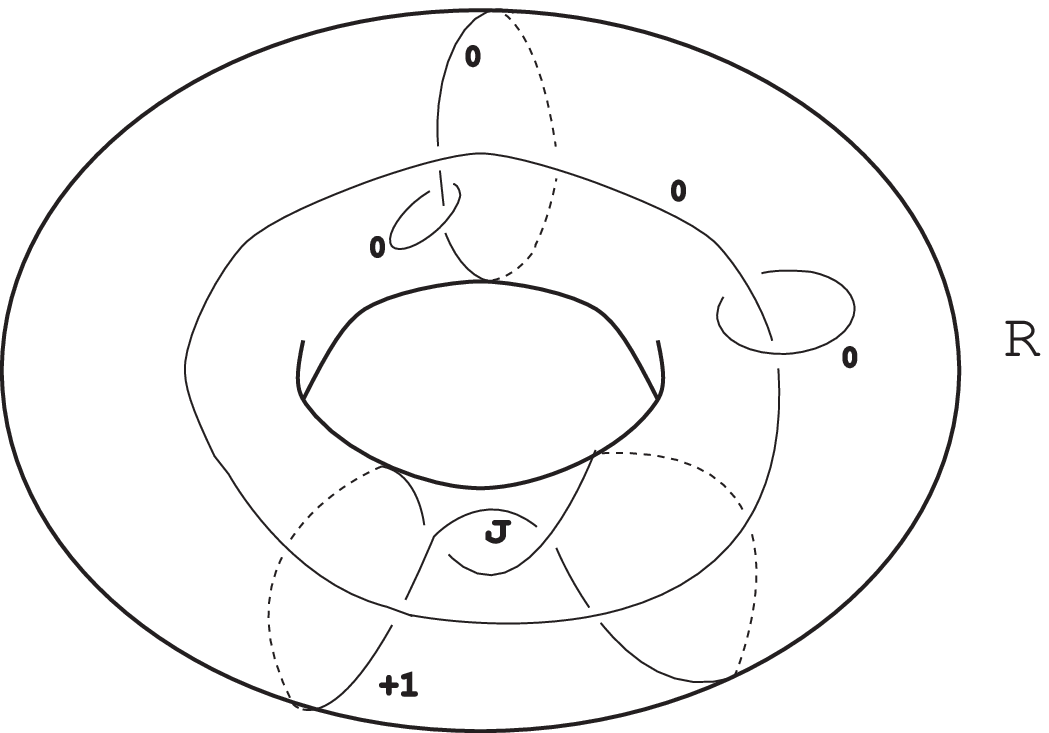}\end{center}
\caption{{\bf 
A framed link in the thickened torus
}\label{figdrawkae}} 
\end{figure}


\begin{figure}
\begin{center}\includegraphics[width=130mm]{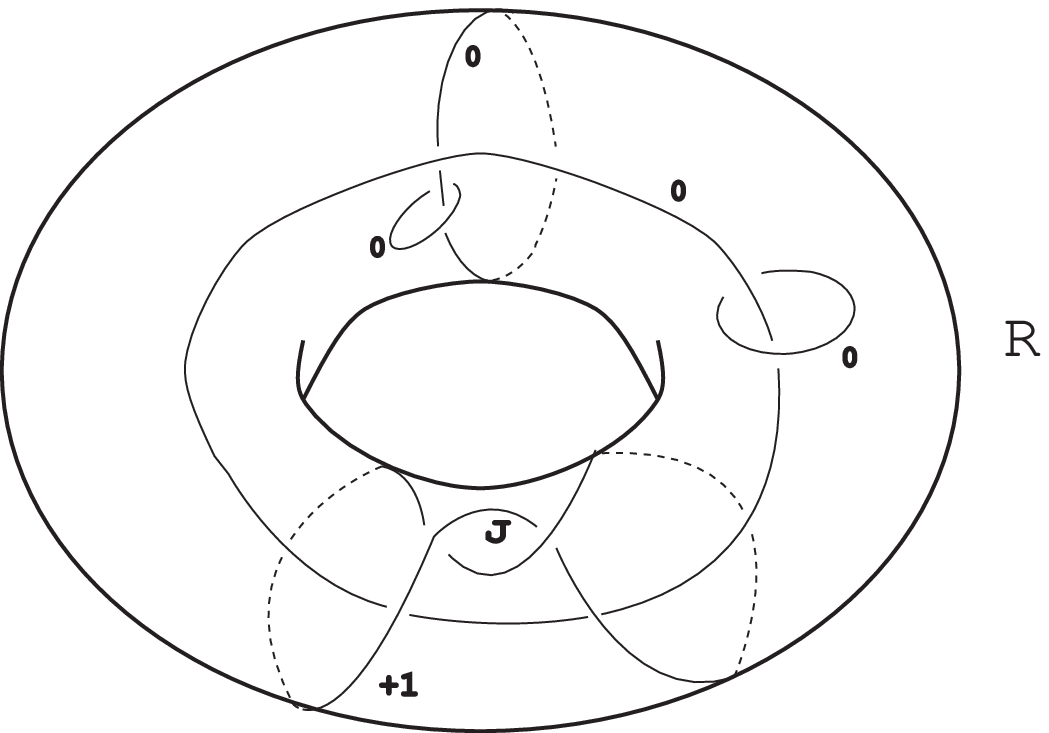}\end{center}
\caption{{\bf 
A framed link in the thickened torus
}\label{figdraw}} 
\end{figure}

\begin{figure}
\includegraphics[width=120mm]{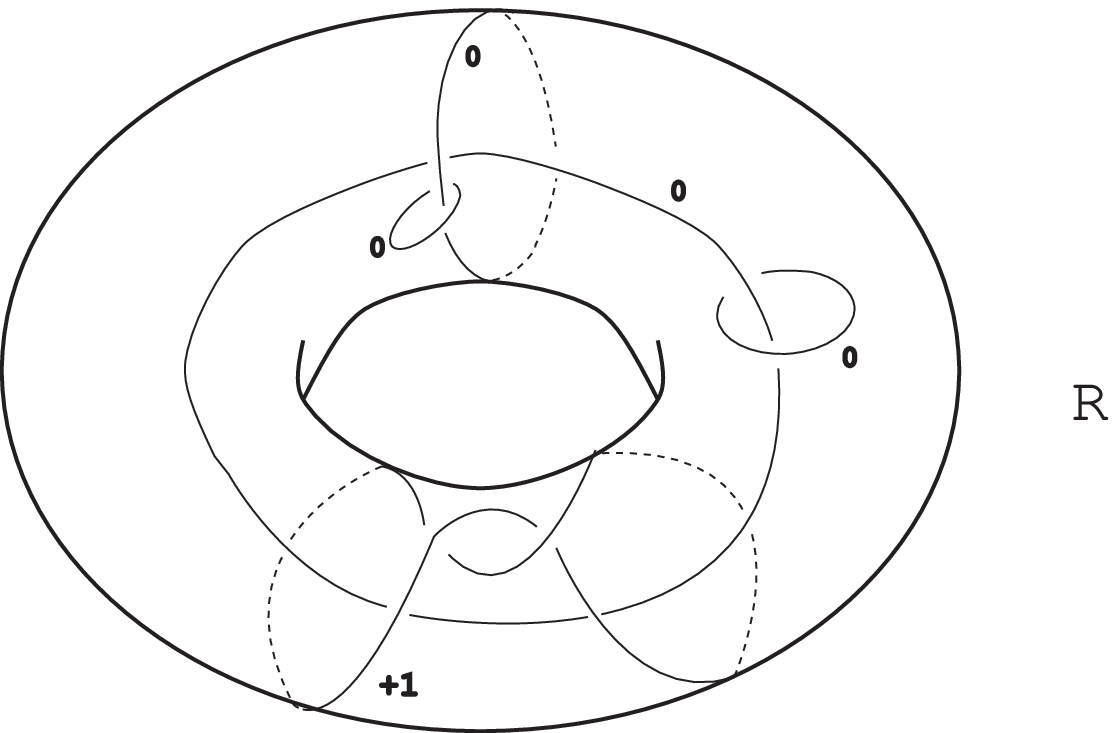}
\caption{{\bf 
A framed link in the thickened torus
}\label{figdrawY}}   
\end{figure}

\begin{figure}
\includegraphics[width=140mm]{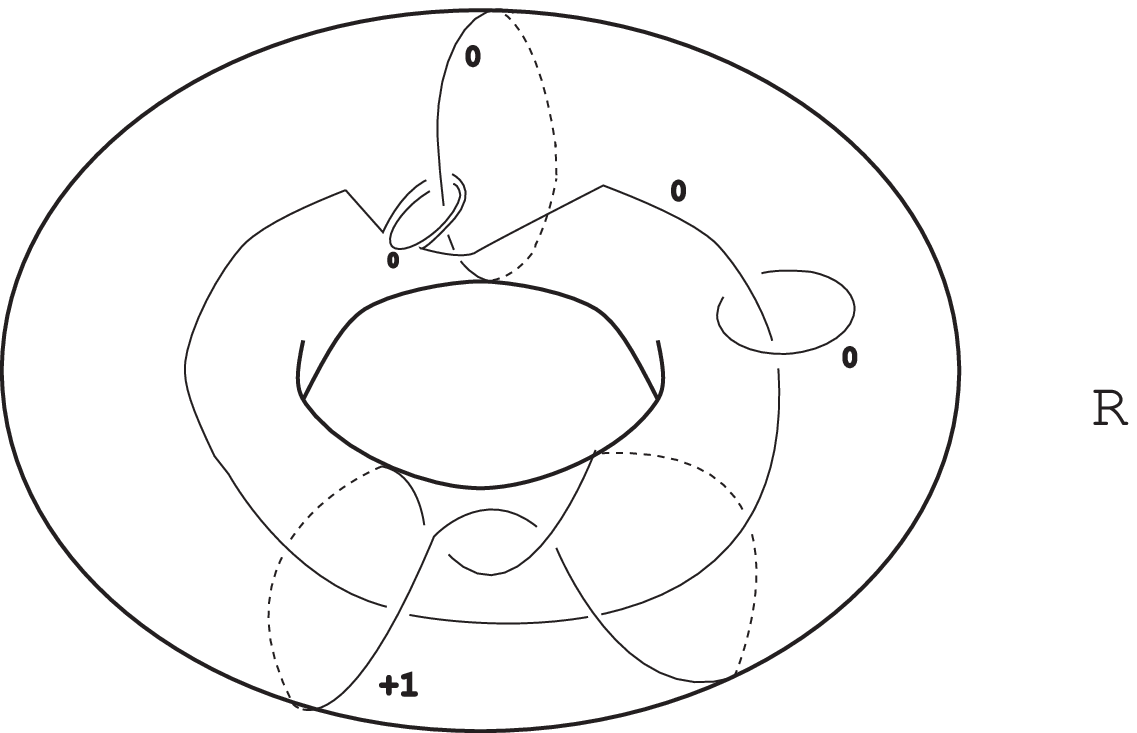}
\caption{{\bf 
A framed link in the thickened torus 
}\label{figdrawX}}   
\end{figure}

\bb
 Figure \ref{fighokuL2} draws a framed link in 
 the complement of  
 the ring-hooked trivial knot, 
 the Hopf link,  
 $(K_1, K_2)$, which is the thickened torus 
 with the standard basis in the boundary.  
Each component has the framing zero.
This framed link represents the thickened torus with the same  symplectic basis condition
because  this framed link is obtained by the $\mathcal O_3$ moves from the empty framed link. 
This framed link satisfies the simple connectivity condition. 
 \\

\begin{figure}
\begin{center}\includegraphics[width=100mm]{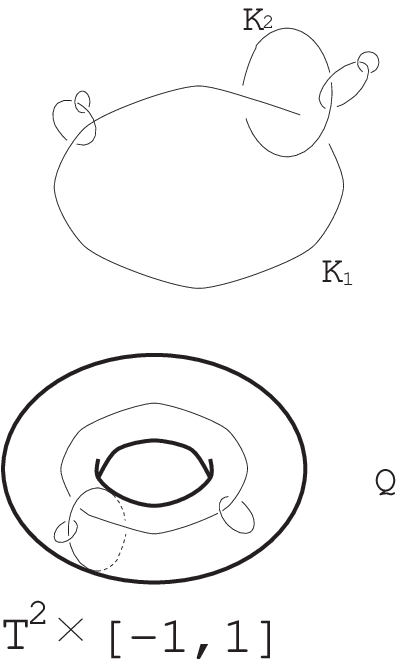}\end{center}
\caption{{\bf 
A framed link in the thickened torus
}\label{fighokuL2}}   
\end{figure}



We prove the following. 
\begin{thm}\label{thmnontri}
\Michi
\end{thm}

\begin{corollary}\label{cor3mfd}
\Mni
\end{corollary}

\h{\bf Proof of Theorem \ref{thmnontri}.} 
Calculate 
our topological quantum knot invariants 
$\upsilon_r$ 
of the right-handed trefoil knot 
by using a framed link drawn in the following figures: 
The lowest  figure  in Figures \ref{fighokuL1}. 
Figures \ref{figdrawkae} - \ref{figdrawX}.  
\\
Calculate 
those of the trivial knot 
by using a framed link drawn in Figure \ref{fighokuL2}. 
They are different by the result of Table \ref{tab:results}. 
\qed
\\

\h{\bf Proof of of Corollary \ref{cor3mfd}.}  
The pair which we calculate in Proof of Theorem \ref{thmnontri} are also
 a pair to prove Corollary \ref{cor3mfd}. 
\qed\\


Therefore 
our topological quantum invariants 
$\upsilon_r$ 
of knots and links are non-trivial invariants. 
Our topological quantum invariants 
$\upsilon_r$ 
of 3-manifolds with boundary are non-trivial invariants. 
\\

We explain how to calculate our invariant 
$\upsilon_r$ 
of 
any other knot $K$ than the trivial knot and the right-handed trefoil knot.  
We construct 
a framed link in the thickened torus with the symplectic basis condition $\mathcal F$, 
which represents the complement of the ring-hooked knot $K$ with the standard basis
as below.
Let  $L$ be any given 2-component link. $L$ may be a ring-hooked knot. 
Take a diagram $D$ in $S^2$ of the Hopf link so that 
 $D$ satisfies the conditions: 
We put the `small trivial' knots with framing  $\pm1$ 
which is null-homologous in $S^3-$(the Hopf link) 
as  in Figure 
   \ref{fighokuK2} 
at some crossing points of $D$.    
Whatever one uses the 3-sphere, 
the thickened torus,  
the fact that each `small knot' is null-homologous, 
this framing is the same (Recall Remark \ref{remreirei}). 
If we carry out surgeries along these `small trivial' knots, 
then we get a diagram of $L$. 

Then the framed link made from all of the above `small trivial' knots 
is in the thickened torus. 
Using the $\mathcal O_3$ moves  
as drawn in Figure \ref{pfigpi1},  
we obtain a framed link which represents 
the complement of the ring-hooked knot with the standard basis 
and which satisfies the simple connectivity condition.   
We make it into a framed virtual link representation in the plane.

\h{\bf Remark:}    
If we regard the Hopf link as a submanifold as drawn in the $S^3$ 
in Figure 
\ref{fighokuK2},  
each `small trivial' knot looks the trivial knot. 
Indeed,  in general, 
the `small trivial' knot is not the trivial knot in the thickened torus.

\section{
\bf 
Our topological quantum invariants 
$\upsilon_r$
is different from  
the Reshetikhin-Turaev topological quantum invariants $\tau_r$
\bf}\label{secchiga}

In  Figure \ref{figchiga1}, 
we draw framed links, 
$L_1$, $L_2$, $L_3$ and $L_4$. 
All components have the framing $0$. 
They satisfy the simple connectivity condition. 
$L_1$ and $L_2$ (respectively, $L_3$ and $L_4$)  
are changed into each other by the $\mathcal O_1$ and $\mathcal O_2$ moves,  
without using a $\mathcal O_3$ move. 
$L_1$ and $L_2$ (respectively, $L_3$ and $L_4$)  
represent the same 3-manifold with boundary.
\\

\h{\bf Remark:}  $Q$ in Figure \ref{fighokuL2} and $L_4$ in Figure \ref{figchiga1} 
are isotopic. 
\\

We prove the following proposition.

\begin{proposition}\label{prKORT} 
In general, each of our quantum invariants 
$\upsilon_r$ 
of 
a 3-manifold with boundary, which is represented by 
a framed virtual link $L_2$ in  Figure \ref{figchiga1}, 
is not equal to that by $L_3$. 
\end{proposition} 

\h{\bf Proof of Proposition \ref{prKORT}.}  Calculation. 
See Table \ref{tab:results}.   
\qed
\\

\h{\bf Remark:} 
Figure \ref{figchiga1} draws 
the framed links, $L_2$ and $L_3$, in the thickened torus. 
In Figure \ref{figchiga1}  the thickened torus  is embedded in the 3-sphere. 
If we forget the thickened torus, 
we make 
framed links, $L_2$ and $L_3$, in the thickened torus in Figure \ref{figchiga1}, 
into 
framed links, $L'_2$ and $L'_3$, in the 3-sphere in Figure \ref{figinS3}, 
respectively. 
Note that the framing of each component of $L_2$ (respectively, $L_3$) 
does not change when we forget the thickened torus  (Recall Remark \ref{remreirei}).
Then 
both framed links $L'_2$ and $L'_3$ 
represent the 3-sphere. 
Therefore they have the same Reshetikhin-Turaev topological quantum invariant. 
Thus our topological quantum invariants 
$\upsilon_r$ 
are different from  
the Reshetikhin-Turaev topological quantum invariants $\tau_r$.
\vs

\begin{figure}
\includegraphics[width=60mm]
{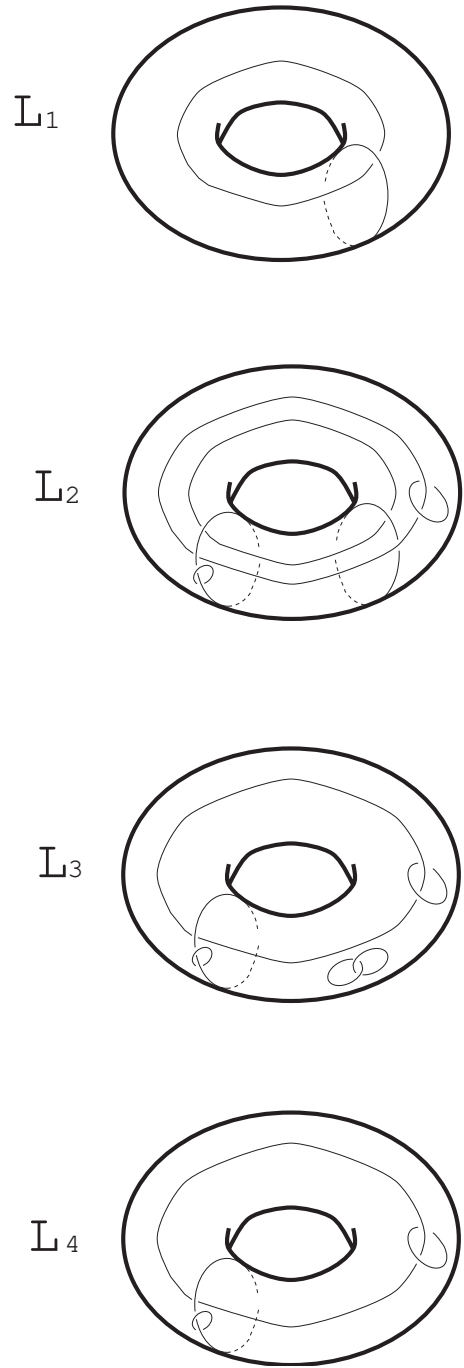}
\caption{{\bf 
Framed links, $L_1$, $L_2$, $L_3$, and $L_4$, in the thickened torus. 
$L_4$ is isotopic to $Q$ in Figure \ref{fighokuL2}.    
%
%
}\label{figchiga1}}   
\end{figure}



\begin{figure}
\includegraphics[width=90mm]{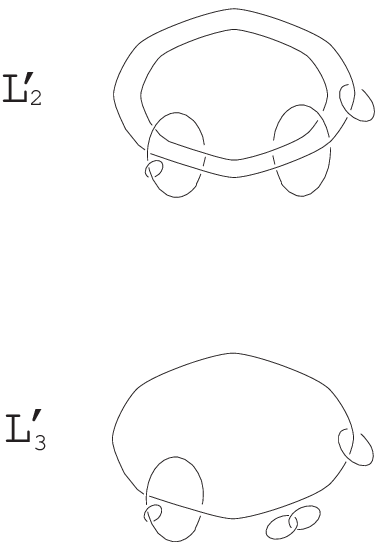}
\caption{{\bf 
Framed links, $L'_2$ and $L'_3$, in the 3-sphere  
}\label{figinS3}}   
\end{figure}

Furthermore, 
there are many ways to embed a thickened surface in $\R^3$. 
Thus it does not determine a topological invariant to use an embedding of thickened surfaces 
in $\R^3$ and 
the Reshetikhin-Turaev topological quantum invariants as above. 

On the other hand, our topological quantum invariants 
$\upsilon_r$ 
are defined 
without using such an embedding.  

\bb

In addition, 
recall  that in the case of our invariants 
$\upsilon_r$ 
of framed links in $F \times I,$ 
the values, 
$\upsilon_r(M)$, 
are invariant under diffeomorphisms of $F \times I.$ 
This is an important property for our invariants 
$\upsilon_r$ 
that does not appear in the usual framework for $WRT$ invariants $\tau_r$. 
 (Remark \ref{remdiff}.)\\

See Figures \ref{figchiga2} and \ref{figchiga3} . 
$A_0$ 
is a framed knot with framing $0$ in the thickened torus. 
If we regard the underlying knot as a virtual knot, it is the virtual right-handed trefoil knot.   

$A_1$ is obtained from $A_0$ by an isotopy.

The framed links, $A_1$ and $A_2$, 
represent the same 3-manifold. 
$A_1$ does not satisfy the simple connectivity condition although $A_2$ satisfies it.

The framed links,  $A_2$ and $J_1$, represent the same 3-manifold. 
Both $A_2$ and $J_1$ satisfies the simple connectivity condition.

$J_2$ is a framed link in the thickened torus with the simple connectivity condition. 
\\

\h{\bf Remark:}  $A_1$ and $A_2$ in Figure \ref{figchiga2} are the same as 
$A_1$ and $A_2$ in Figure \ref{figcheck}, respectively. \\

We prove the following proposition.

\begin{proposition}\label{proKORT2}
Let $J_1$ and $J_2$ be as shown in Figure $\ref{figchiga3}$ 
and as discussed above.
In general, each of  our quantum invariants 
$\upsilon_r$ 
of a 3-manifold with boundary, 
which is represented by $J_1$,  
is not equivalent to that by $J_2$. 
\end{proposition}

\h{\bf Proof of Proposition \ref{proKORT2}.}
Calculation. Table \ref{tab:results}.   
\qed\\

\h{\bf Remark:}  
Figure \ref{figchiga3} draws 
the framed links, $J_1$ and $J_2$. in the thickened torus. 
In Figure \ref{figchiga3}  the thickened torus  is embedded in the 3-sphere. 
If we forget the thickened torus, 
we make 
framed links, $J_1$ and $J_2$, in the thickened torus in Figure \ref{figchiga3}, 
into 
framed links, $J'_1$ and $J'_2$, in the 3-sphere in Figure \ref{figinS3p}, 
respectively. 
Note that the framing of each component of $J_1$ (respectively, $J_2$) 
does not change when we forget the thickened torus  (Recall Remark \ref{remreirei}).
Then 
both the framed links $J'_1$ and $J'_2$ 
represent the same 3-manifold. 
Therefore they have the same Reshetikhin-Turaev topological quantum invariant. 
Thus our topological quantum invariants 
$\upsilon_r$
are different from the Reshetikhin-Turaev topological quantum invariants $\tau_r$.
\vs


\begin{figure}
\includegraphics[width=90mm]{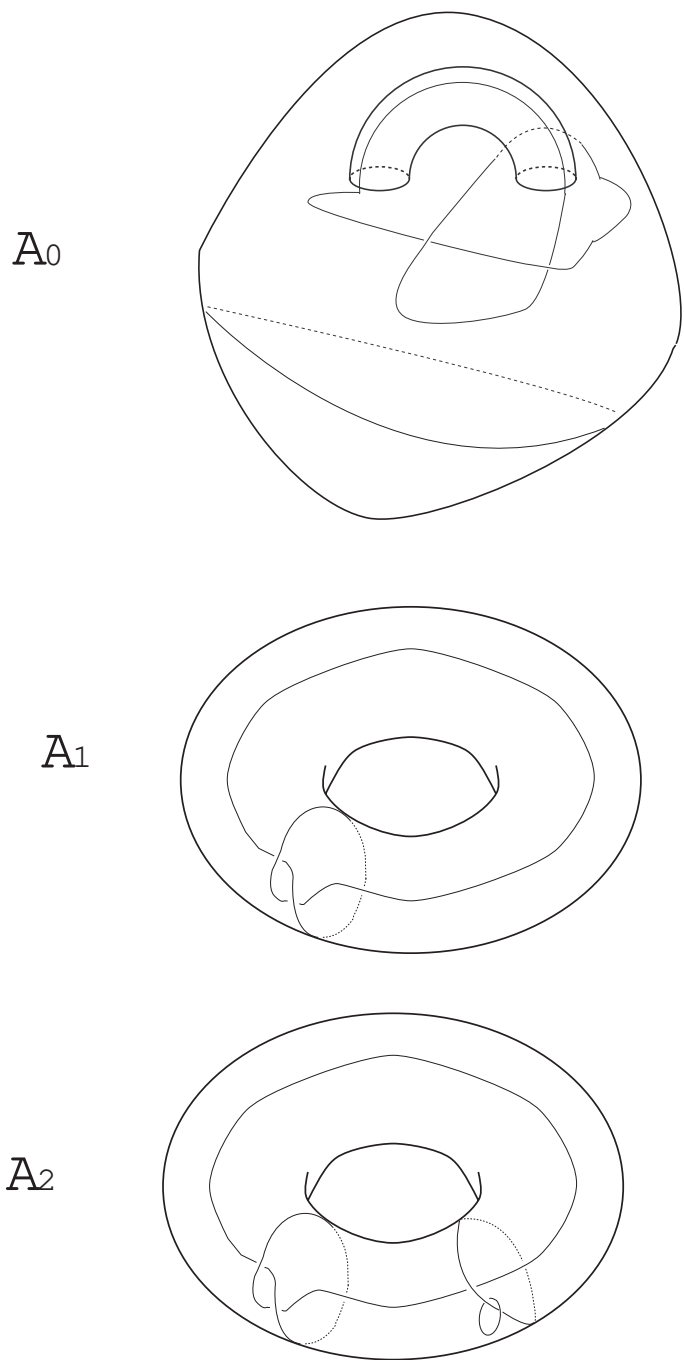}
\caption{{\bf 
Framed links in the thickened torus 
}\label{figchiga2}}   
\end{figure}



\begin{figure}
\includegraphics[width=100mm]{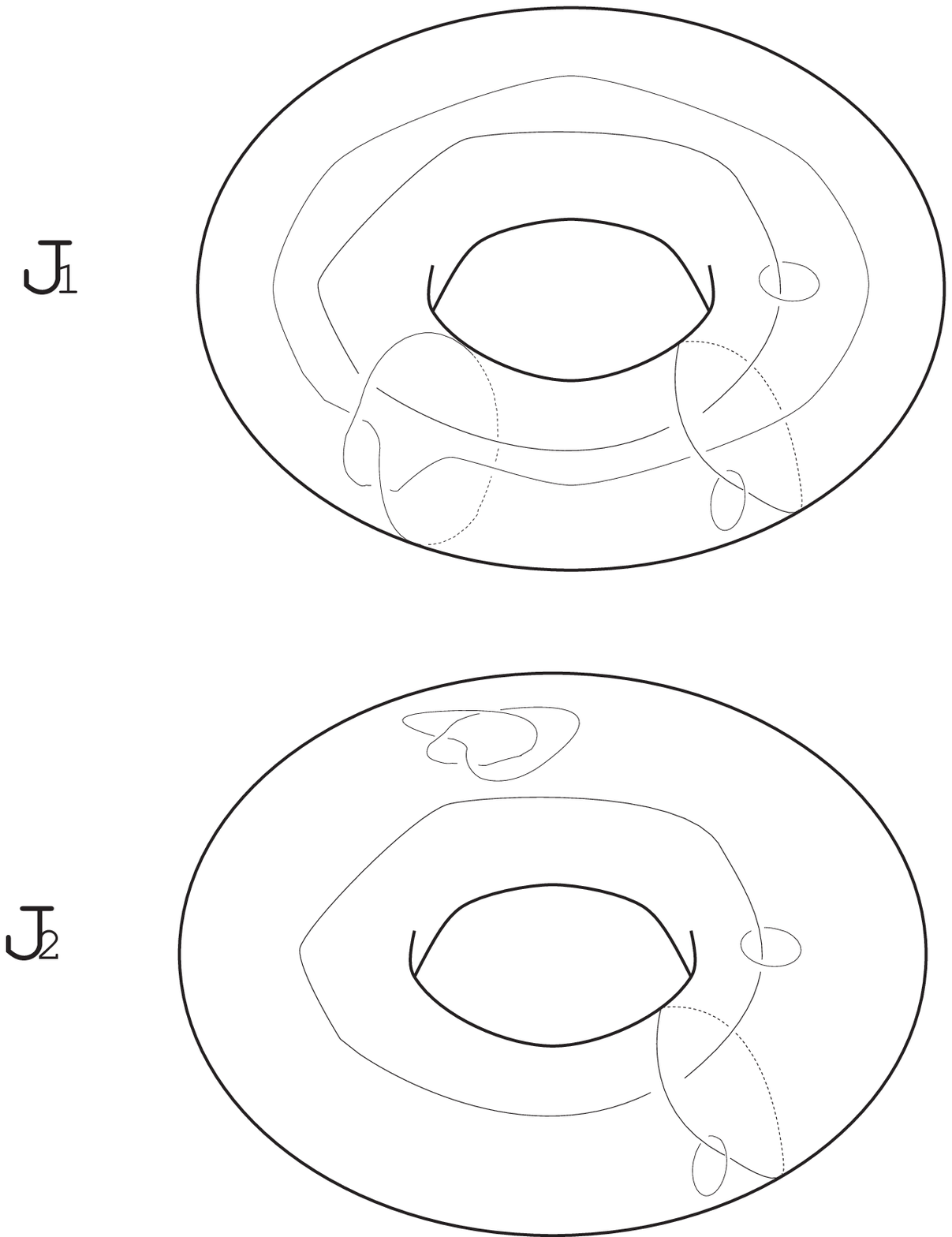}
\caption{{\bf 
Framed links in the thickened torus 
}\label{figchiga3}}   
\end{figure}

\begin{figure}
\includegraphics[width=100mm]{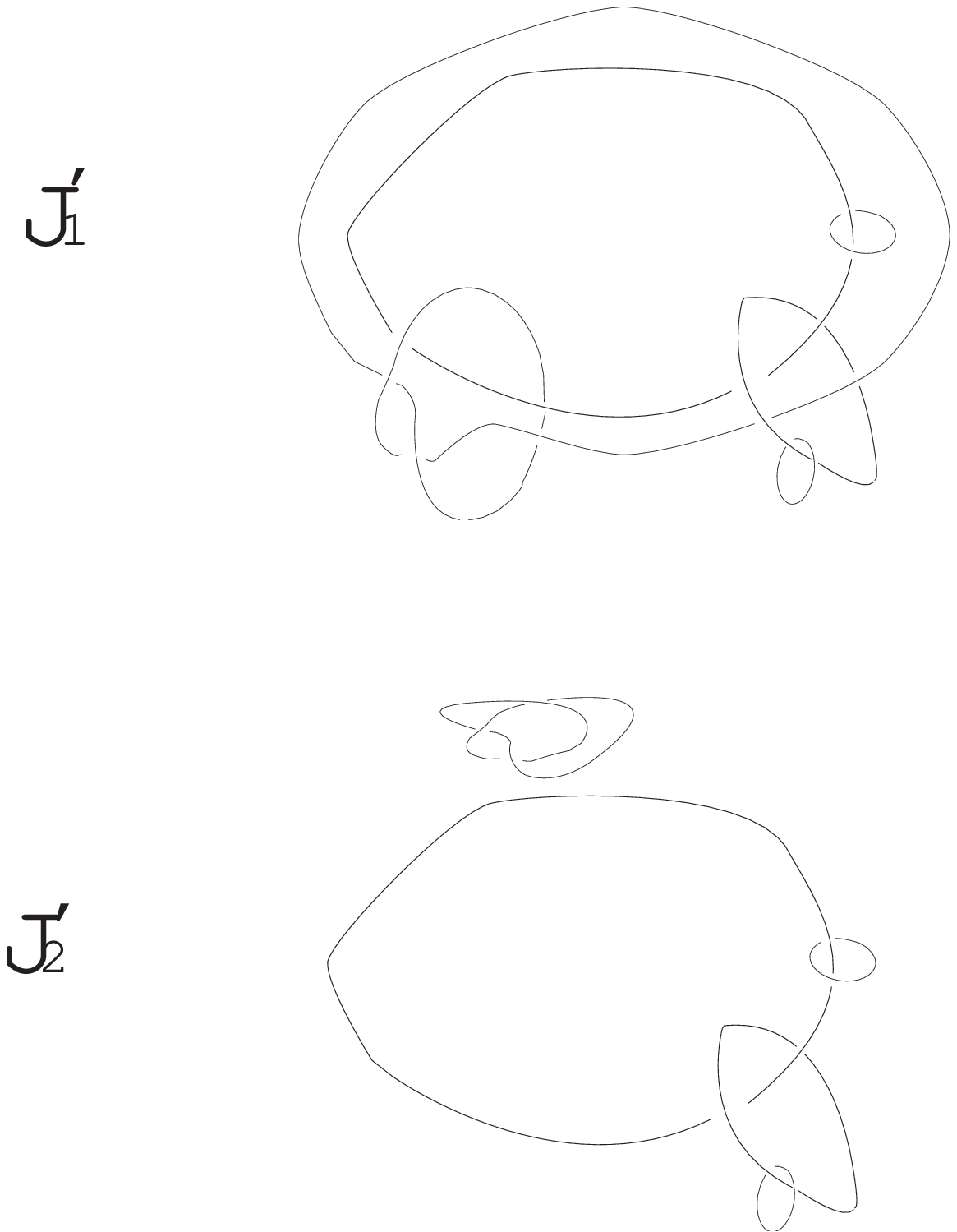}
\caption{{\bf 
Framed links, $J'_1$ and $J'_2$, in the 3-sphere  
}\label{figinS3p}}   
\end{figure}









\section{\bf 
Calculation result}\label{seckekka} 

We calculate examples which we exhibited in this part, Part \ref{part2},  explicitly. 
See Table \ref{tab:results}. \\

The closure of the q-symmetrizer, denoted as $\Delta_n$ is the closure of a linear combination of 
elements in the Temperly-Lieb algebra and 
\begin{equation}
\Delta_n = (-1)^n \frac{A^{2n} - A^{-2n}}{A^2 - A^{-2}}.
\end{equation}
The recursive definition is utilized to evaluate the Kauffman-Ogasa invariant. 
The values of $\Delta_n$ for $n=0$ and $n=-1$ are formally defined:  $ \Delta_{-1} = 0$ and $ \Delta_0 =1$. The symbol $\Delta_1 $ represents a single closed loop and $ \Delta_1 = 1$. Then for $n>0$, 
\begin{equation}
\Delta_{n+1} = d \Delta_n - \Delta_{n-1}. 
\end{equation}

Each crossing in the knot diagram is transformed into a pair of trivalent vertices with labeled edges; shown in figure \ref{fig:crossingtovertices}. 
 
\begin{figure}
\[ \begin{array}{c} \includegraphics[scale=0.35]{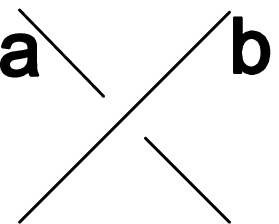} \end{array}
= \sum_{Adm(a,b,i)} \frac{ \Delta_i}{\theta (a, b, i) } \lambda_i ^{(ab)  } 
\begin{array}{c} \includegraphics[scale=0.3]{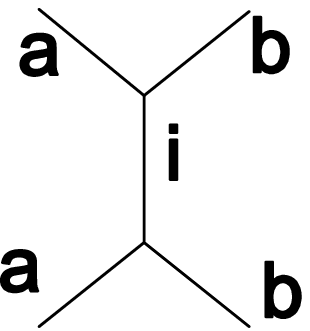} \end{array}
\]
\caption{Converting link diagrams into graphs}
\label{fig:crossingtovertices}
\end{figure}

The labeled trivalent vertex represents a tangle as shown using the q-symmetrizer in figure \ref{fig:trivalentvertex}. 
\begin{figure}
\[ \begin{array}{c} \includegraphics[scale=0.3]{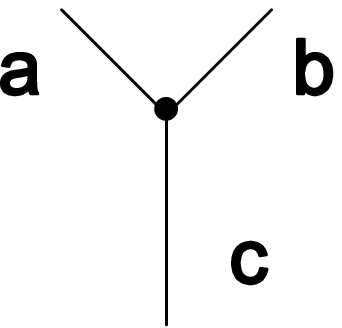} \end{array} \leftrightarrow
\begin{array}{c} \includegraphics[scale=0.3]{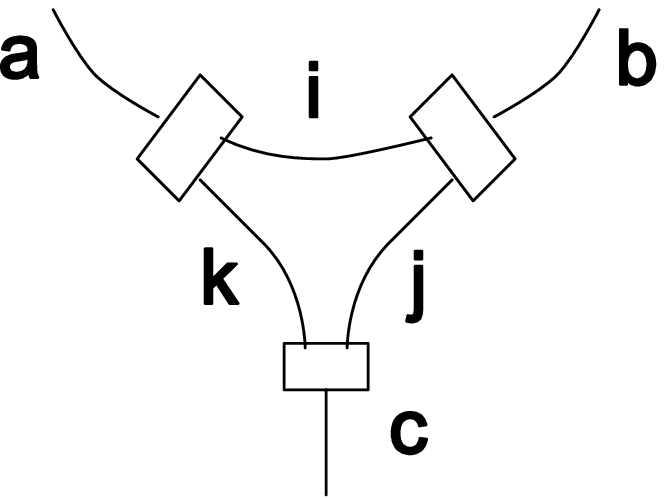} \end{array} 
\]
\caption{Expanding the trivalent vertex}
\label{fig:trivalentvertex}
\end{figure}
The edge labels must satisfy the equations
\begin{align}\label{vertexadmiss} 
i&= (a+b-c)/2, & j&= (a+c -b)/2, & k&=(b+c -a)/2 .
\end{align}

The labeled link diagram is converted into a a sum of weighted trivalent graphs where each labeled edge represents a q-symmetrizer. The graph can be successively simplified using the theta net and tetrahedral formulas given early in the paper.  This is
evaluated to a complex number by fixing an integer $r > 2$ and letting $A=e^{\frac{\pi i }{2r}}$.  In this case,   $\Delta_{r-1}=0$.   
This, combined with the restriction on trivalent vertices, leads admissibility restrictions on the labels of the diagrams.  Non-admissible labels evaluate to zero. 

We summarize the admissibility requirement for any vertex in a diagram. 
For a fixed integer $r>2$, any labeled edge must satisfy $0 \leq n \leq r-2$. 
Combined with the vertex admissibility condition in equation \ref{vertexadmiss}, we obtain that
$a+b+c \leq 2r-4$ and $i, j, k \geq 0 $.

For a virtual link diagram $K$ obtained from a link in $F \times I$ and subject to the simple boundary condition, let $c$ denote the number of components is $c$ and let $B$ denote the linking matrix.  We let $b_-$ (respectively $b_+$) denote the number of negative (respectively positive) eigenvalues of $B$. Then
\begin{equation}
Z(K) = \langle {\omega * K} \rangle \mu^{c+1} \alpha ^{-(b_+ - b_-)}
\end{equation}
where
\begin{align*}
\mu &= \sqrt{ \frac{2}{r} sin} \left( \frac{\pi}{r} \right), & \alpha &= (-i)^{r-2} e^{i \pi ( \frac{3 (r-2)}{4r}}.
\end{align*}
For a single component diagram $K$, 
$\langle \omega * K \rangle = \sum_{i=0} ^{r-2} \Delta_i K^i$ where $K$ is labeled with $i$. 
Using this normalization, the $Z(U) = 1$ where $ U$ is the zero framed unknot.

Working with virtual link diagrams, we obtain elements of the virtual Temperly-Lieb algebra instead of the Temperly-Lieb algebra. As a result, our simplified diagrams contain virtual crossings and we are sometimes unable to completely reduce the diagrams to single loops.   In the case that a reduced diagram has virtual crossings, we use symmetry and recursive evaluations of the q-symmetrizer to obtain a complex value. 

Consider the virtual Hopf link $H$. 
Let $a$ and $b$ denote the labels placed on the components of the link. The linking matrix is
\begin{equation*}
B=\begin{bmatrix} 0 & -1 \\ 0 & 0 \end{bmatrix}.
\end{equation*} 
Then
\begin{equation}
Z(H) = \langle \omega * H \rangle \mu^3 \alpha^0.
\end{equation}
We expand $\langle \omega * H \rangle$:
\begin{equation} 
\langle \omega * H \rangle = \sum_{a,b=0} ^{r-2} \Delta_a \Delta_b  
\begin{array}{c} \includegraphics[scale=0.35]{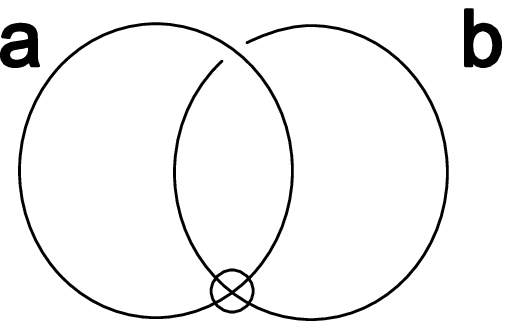}\end{array} .
\end{equation}
We then exchange the crossing for an edge, constructing a trivalent graph
\begin{equation} \label{equation:hopfred}
\sum_{Adm(a,b,i)} \Delta_a \Delta_b \frac{ \Delta_i}{\theta (a,b,i)} \lambda_i ^{ab} 
\begin{array}{c} \includegraphics[scale=0.35]{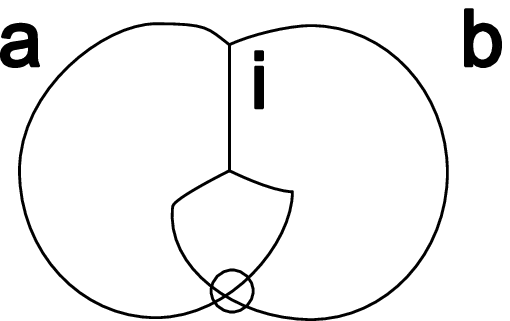} \end{array} .
\end{equation}

This diagram cannot be simplified  further using the existing reduction formulas. 
The graph represents a collection of closed curves with virtual crossings. The graph in 
\ref{equation:hopfred} is symmetric. Under the detour move, any pair of edges in the graph can contain the virtual crossing. We refer to this diagram as a \textit{twisted theta}, $\tilde{ \theta} (a,b,i)$. 

Finally, we obtain the formula
\begin{equation}
\langle \omega * H \rangle = \sum_{Adm(a,b,i)} \Delta_a \Delta_b  \Delta_i \frac{ \tilde{\theta } (a,b,i) }{\theta (a,b,i)} \lambda_i ^{ab}  .
\end{equation}

Note that $\tilde{ \theta}(a,a,0) $ evaluates to $ \Delta_a$ for all values of $a$. Other values must be evaluated by hand. 

We consider another example that cannot be reduced using the existing reduction formulas. 

We evaluate the diagram $A_2$ from figure \ref{figchiga2}. This three component link has a linking matrix of the 
form:
\begin{equation}
\begin{bmatrix}
2 & 0 & 0 \\ 1 & 0 & 1 \\ 0 & 1 & 0 \end{bmatrix}.
\end{equation} 
Note that there are 2 positive eigenvalues and 1 negative eigenvalue. 
Therefore, 
\begin{equation}
Z(A_2) =  \mu^3 \alpha^{-1} \langle \omega * A2 \rangle. 
\end{equation} 
We reduce and evaluate $\langle \omega * A_2 \rangle$.
\begin{equation}
\langle \omega * A_2  \rangle = \sum_{a,b,c=0} ^{r-2} \Delta_a \Delta_b \Delta_c \left\langle \begin{array}{c} \includegraphics[scale=0.35]{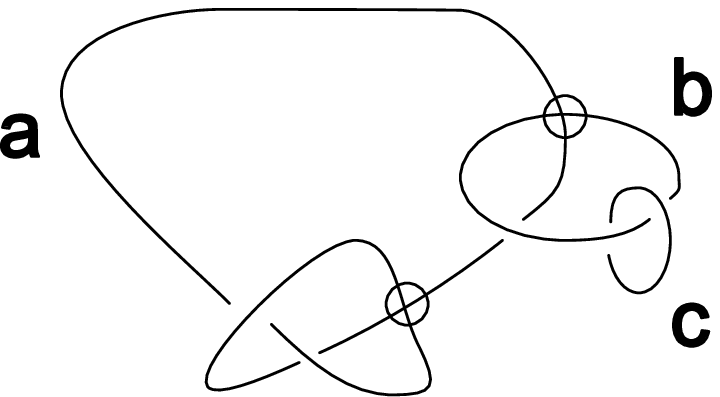} \end{array} \right\rangle.
\end{equation}

Next, }

\begin{gather}
\langle \omega * A2  \rangle = \sum_{
\begin{matrix} Adm(a,a,j) \\  Adm(a,a,m) \\ Adm(b,c,i) \end{matrix}}  \Delta_a  \Delta_c
\Delta_i (\lambda_i ^{bc})^2 (\lambda_j ^{aa})^2 (\lambda_n ^{aa})^2 \frac{ \Delta_m ^2}{\theta(a,a,m)^3} \frac{\Delta_j}{ \theta(a,a,j)} \\ \nonumber
  Tet \begin{bmatrix} a & a & m \\ a & a & j \end{bmatrix}  
 \left\langle \begin{array}{c} \includegraphics[scale=0.35]{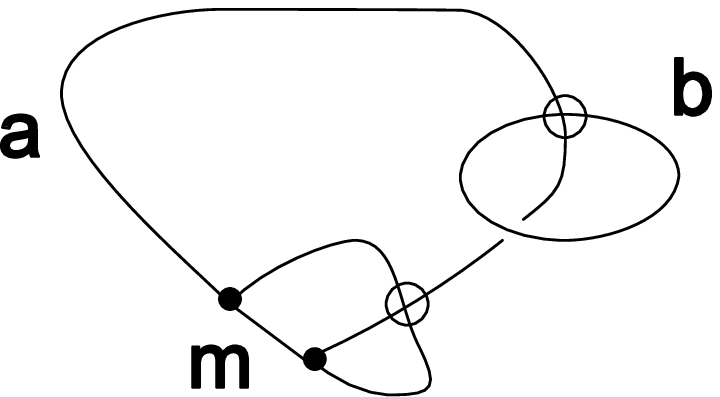}  \end{array} \right\rangle.
\end{gather}

{
The final reduction to a spatial graph results in the formula:}
\begin{gather*}
\langle \omega * A2  \rangle = \sum_{\begin{matrix} Adm(a,a,j) \\ Adm(a,a,m) \\  Adm(b,c,i) \\ Adm(a,b,p) \end{matrix} } \Delta_a  \Delta_c
\Delta_i (\lambda_i ^{bc})^2 (\lambda_j ^{aa})^2 (\lambda_n ^{aa})^2 \frac{ \Delta_m ^2}{ \theta(a,a,m)^3}  \frac{\Delta_j} { \theta(a,a,j)}  \frac{ \Delta_p}{\theta(a,b,p)} \\ \nonumber
 Tet \begin{bmatrix} a & a & m \\ a & a & j \end{bmatrix}  
 \left\langle \begin{array}{c} \includegraphics[scale=0.35]{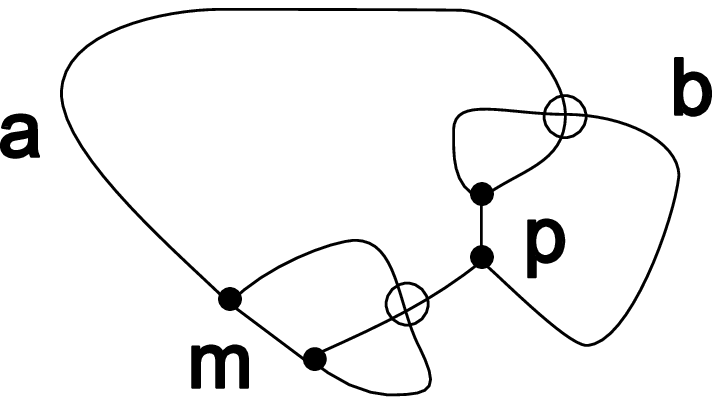} \end{array} 
 \right\rangle.
\end{gather*}

{
Note that the reduced graph has 4 vertices and forms a tetrahedral net  structure with virtual crossings.  This graph must be evaluated by hand to obtain an element of the algebra $ \mathbb{Z} [A, A^{-1}]$. Fortunately, the figure has some symmetries, simplifying the evaluation.

Our computational results are in Table \ref{tab:results}.

}
{
Future questions for research include finding a closed formula for $\tilde{\theta} (a,b,i)$ and the virtual tetrahedrons.

\bb
We calculate our topological quantum invariants 
$\upsilon_r$ 
of the following 3-manifolds with the boudary condition $\mathcal B$.
They are represented by framed links in the thickened torus 
with the symplectic basis condition $\mathcal F$.
\\

$A_1, A_2$ in Figure \ref{figcheck}. 
\\

$Y, C$ in 
Figure \ref{figcheck00Y}.\\

$R$ in Figures 
\ref{fighokuL1} 
and 
 \ref{figdrawkae}
- \ref{figdrawX}.
\\

$Q$ in Figure \ref{fighokuL2}.
\\

$L_1, L_2, L_3, L_4$  in Figure \ref{figchiga1}. 

$\upsilon_r(L_1)=\upsilon_r(L_2)$. 

$\upsilon_r(L_3)=\upsilon_r(L_4)$. 

$L_4$ and $Q$ are isotopic. Hence 
$\upsilon_r(L_4)=\upsilon_r(Q)$. 
\\

$A_0, A_1, A_2$ in 
Figure \ref{figchiga2}. 

$A_0$ and $A_1$ are isotopic. Hence 
$\upsilon_r(A_0)=\upsilon_r(A_1)$. 

$A_1$ and $A_2$ in Figure \ref{figchiga2} are the same as $A_1$ and $A_2$ in Figure \ref{figcheck}, respectively. 
\\

$J_1, J_2$ in Figure \ref{figchiga3}. 

$\upsilon_r(A_2)=\upsilon_r(J_1)$.
\\ 

$X$, $A$, and $B$ are introduced right below.
\\

\bb
{\normalsize 
\begin{table}[h!]
  \begin{center}
    \begin{tabular}{l|c|c|c|c} 
      \textbf{Framed link} 
       & \textbf{$r=3$} & \textbf{$r=4$} & \textbf{$r=5$} \\
 
      \hline
 
      
      $L_1$, $L_2$ 
      &1.06066 - 0.353553 $i$ & 0.0967185 + 1.20711 $i$ & 0.553238 + 1.04288 $i$ \\
   
   
      
      $L_3, L_4, Q, X, Y, A, B, C$ 
      &  0.707107 & 0.5 & 0.37148 \\

       $A_1$ 
       &  0.707107 + 0.707107 $i$ & -1.30656 + 0.92388 $i$ & 
      -1.58479 - 1.72679 $ i$  \\


       $A_2$, $J_1$ 
       & -0.25 + 0.103553 $i$ &0.0544203 - 0.253256 $i$ 
      & 
      \\

      $J_2$  
      & 0.707107 $i$ & 0.353553 - 0.353553 $i$ & -0.158114 + 0.716377 $i$\\

$R$ 
    &  -0.707107  
    & 0
     &  0.352125 - 0.484658i
      \\
      
    \end{tabular}
\bigbreak   \center \caption{Computational Results.  
}
  \label{tab:results}
    \end{center}
\end{table}
}

We ask a question. 
See Figure \ref{figcheck00}. 
\begin{figure}
\includegraphics[width=160mm]{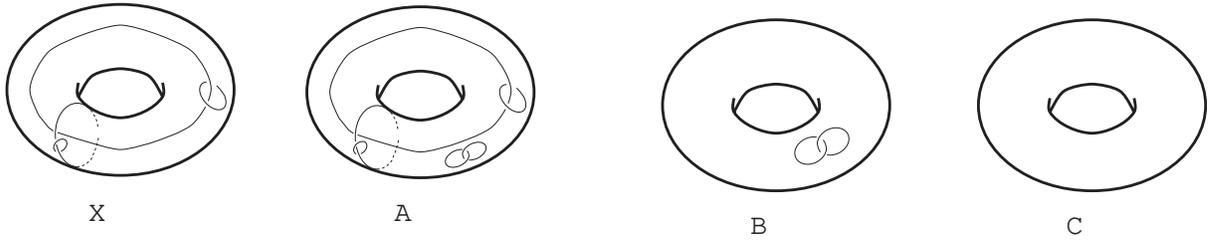}
\caption{{\bf  
Framed links in thickened tori. All framings are zero. 
Note that $A$ in Figure \ref{figcheck00} is $L_3$ in Figure \ref{figchiga1}, 
and that $X$ in Figure \ref{figcheck00}  is $L_4$ in Figure \ref{figchiga1}.
}\label{figcheck00}}   
\end{figure}
$C$ represents the empty framed link. 
$X$, $A$, $B$, and $C$ represent the same 3-manifold. \\
$X$ is obtained form $A$ by a sequence of only the $\mathcal O_1$ and the $\mathcal O_2$ moves. \\
$B$ is obtained from $C$  by a sequence of only the $\mathcal O_1$ and the $\mathcal O_2$ moves. \\ 
$A$ is not obtained from $B$ without using the $\mathcal O_3$ move 
because their fundamental groups are different. 
Therefore $A$ is not obtained form $B$ by 
a sequence of the $\mathcal O_1$ and the $\mathcal O_2$ moves.  
\\ 
$X$ and $A$ satisfy the simple connectivity condition while $B$ and $C$ do not.

Each of the Dye-Kauffman quantum invariants 
$\varsigma_r$ 
of $X$ and 
that of $A$ are the same because $X$ is obtained from $A$ 
by the $\mathcal O_1$ and the $\mathcal O_2$ moves. 

Each of the Dye-Kauffman quantum invariants 
$\varsigma_r$ 
of $B$ and 
that of $C$ are the same because $B$ is obtained from $C$ 
by the $\mathcal O_1$ and the $\mathcal O_2$ moves. 
\\


$C$ in Figure \ref{figcheck00} is the same as $C$ in Figure \ref{figcheck00Y}. 

$\upsilon_r(X)=\upsilon_r(A)$. 

$\upsilon_r(B)=\upsilon_r(C)$.

 $A$ in Figure \ref{figcheck00} is $L_3$ in Figure \ref{figchiga1}. 
 
$X$ in Figure \ref{figcheck00} is $L_4$ in Figure \ref{figchiga1}. 
\\

In the case of $r=3,4,5$, 
the Dye-Kauffman invariants 
$\varsigma_r$ 
of $X,$ $A$, $B$, and $C$ 
have the same apparent value.   
See Table \ref{tab:results}.   
\\

\h
{\it Question}. 
Are the Dye-Kauffman quantum invariants 
$\varsigma_r$ 
of $X$,  $A$, $B$, and $C$ the same
for all $r$?

\vs

\section*{\bf 
Appendix}\label{secApp}

 We discuss some open questions.

 Reshetikhin and Turaev 
\cite[Theorem 3.3.3, page 560]{RT}
defined  invariants for links in a closed oriented 
3-manifold $M$. These invariants depend on the use of the colored Jones polynomials at roots of unity and the use of the Kirby calculus.
If $M=S^3$, these invariants are invariants of links in $S^3$.

Apply the definition in \cite[Theorem 3.3.3, page 560]{RT}
strictly to the case of $S^3$. 
The readers can understand easily 
that the definition of these invariants 
is different from that of the Jones polynomial. 

It is an open question whether 
these invaraints retrieves the Jones polynomial. 

The following two questions both are open.
\\

\h{\bf Qustion A.1.}
Suppose that these invariants of a link $\mathcal L$ in $S^3$ 
and those of $\mathcal L'$ in $S^3$  
are the same. (We may not be able to know the coincidence by a finite times of operation. 
We may just suppose this condition abstractly.)  
Then are the Jones polynomial of $\mathcal L$ and that of $\mathcal L'$ the same?
\\

\h{\bf Qustion A.2.}
By a finite times of 
explicit calculation of  (a partial information of) these invariants
of a link $\mathcal L$ in $S^3$,    
can we determine the Jones polynomial of  $\mathcal L$? 
\\

We have not known whether  
 the Reshetikhin-Turaev invariants of links in any 3-manifold is 
an extension of the Jones polynomial of links in $S^3$.

\bb
We have another open question.
\\

\h{\bf Qustion A.3.}
Calculate Witten's well-known path integral 
for links in other 3-manifolds than $S^3$. 
\\

\h{\bf Remark:} 
Witten calculated only the $S^3$ case. See \cite{W}.

\bb

Even in the `physics' level,  the Jones polynomial has not been extended 
to all 3-manifold cases.
Witten only wrote a Lagrangian and an observable for a path integral, 
in the case of links in all (closed oriented) 3-manifolds. 
Only writing a path integral 
never means that the path integral has been calculated. 
Before calculating it explicitly (in the `physics' level), 
we cannot say that the path integral defines a value. 

Recall a current situation of QCD and a history of QED. 
Before Tomonaga, Feynman and Schwinger discovered renormalization, 
we wrote a well-known Lagrangian and wrote path integral for QED. 
We  write a well-known Lagrangian and write path integral for QCD, 
but we cannot say that we complete QCD. 

\bb

It is also an open question whether 
Khovanov homology (\cite{B, Khovanov1999})
and 
 Khovanov-Lipshitz-Sarkar homotopy type (\cite{LSk, LSs, LSr, Seed})
are extended to 
the case of links in any other 3-manifold than the 3-sphere. 
Both are now extended only  
to the case of thickened surfaces. 
See the case of Khovanov homoology for 
Asaeda, Przytycki, and Sikora \cite{APS}, 
Manturov \cite{Man} (arXiv 2006), 
Rushworth \cite{Ru},   
Tubbenhauer \cite{Tub},  
and  
Viro \cite{Viro}.
Manturov and Nikonov \cite{MN} refined the definition of \cite{APS}.  
Dye, Kastner, and Kauffman \cite{DKK},   
Nikonov \cite{Igor}, 
and 
Kauffman and Ogasa \cite{KauffmanOgasasq} refined the definition of \cite{Man}.
Kauffman and Ogasa 
\cite{KauffmanOgasasq} extended 
the second Steenrod square acting on Khovanov homology for links in $S^3$ to the case of thickened surfaces. 
Kauffman, Nikonov and Ogasa 
\cite{KauffmanNikonovOgasa, KauffmanNikonovOgasaT2} 
extended Khovanov-Lipshitz-Sarkar stable homotopy type to the case of thickened surfaces.  
Some of them are defined by using virtual links.

\vs


\vs

\noindent
Heather A. Dye

\noindent
Division of Science and Math

\noindent
McKendree University

\noindent
701 College Rd.

\noindent
Lebanon, IL 62254

\noindent
USA

\noindent
heatheranndye@gmail.com
\\

\noindent
Louis H. Kauffman

\noindent
Department of Mathematics, Statistics and Computer Science

\noindent
University of Illinois at Chicago

\noindent
851 South Morgan Street

\noindent
Chicago, Illinois 60607-7045

\noindent
USA

\noindent
kauffman@uic.edu
\\

\noindent
Eiji Ogasa

\noindent
Meijigakuin University, Computer Science 
 
\noindent
Yokohama, Kanagawa, 244-8539 

\noindent
Japan 

\noindent
pqr100pqr100@yahoo.co.jp  

\noindent
ogasa@mail1.meijigkakuin.ac.jp

}

\begin{thebibliography}{ABCD}
{



\bibitem{APS} 
M. M. Asaeda, J. H. Przytycki, and A. S. Sikora; 
Categorification of the Kauffman bracket skein
module of $I$-bundles over surfaces,  
{\it Algebraic \& Geometric Topology}
4 (2004) 1177–1210 ATG


\bibitem{B}
D. Bar-Natan: 
On Khovanov's categorification of the Jones polynomial, 
{\it Algebr. Geom. Topol.} 2(2002), 337–370 (electronic). MR 1917056 (2003h:57014). 





\bibitem{DKK} 
H A Dye, A Kaestner, and L H Kauffman:  
Khovanov Homology, Lee Homology and a Rasmussen Invariant for Virtual Knots, 
{\it Journal of Knot Theory and Its Ramifications} 26 (2017).

\bibitem{DK}
H. A. Dye, and L. H. Kauffman: 
Virtual Knot Diagrams and the Witten-Reshetikhin-Turaev Invariant, 
 arXiv:math/0407407
 
\bibitem{FR} 
R. A. Fenn, and C. P. Rourke: 
On Kirby’s calculus of links, 
{\it Topology} 18 (1979), 1-15.
 
 
 
 
\bibitem{IMvbunrui}
D. P. Ilyutko and V. O. Manturov:   
Virtual Knots:
The State of the Art, 
{\it World Scientific Publishing Co. Pte. Ltd.} 2012. 
English translation of \cite{Mvbunrui}.


 
\bibitem{Jones} 
 V. F. R. Jones: Hecke Algebra representations of braid groups and link   
{\it Ann. of Math.} 126 (1987) 335-388. 



\bibitem{Kauffman1} 
L. H. Kauffman: 
Talks at MSRI Meeting in January 1997, AMS Meeting at University of Maryland, College Park in March 1997, Isaac Newton Institute Lecture in November 1997, Knots in Hellas Meeting in Delphi, Greece in July 1998, APCTP-NANKAI Symposium on Yang-Baxter Systems, Non-Linear Models and Applications at Seoul, Korea in October 1998


\bibitem{knotphys}
L. H. Kauffman: 
Knots and Physics
{\it Series on Knots and Everything, Vol. 1, World Scientific} 1991, 1994, 2001.

\bibitem{Kauffman} 
L. H. Kauffman: 
Virtual Knot Theory, 
{\it Europ. J. Combinatorics} (1999) 20, 663–691, 
{\it Article No. eujc}.1999.0314, 
{\it Available online at} http://www.idealibrary.com 
math/9811028 [math.GT].





\bibitem{Kauffmani} 
L. H. Kauffman: 
Introduction to virtual knot theory 
{\it J. Knot Theory Ramifications} 21 (2012), no. 13, 1240007, 37 pp.


\bibitem{tl}
L. H. Kauffman and S. L. Lins: 
Temperly-Lieb Recoupling Theory and Invariants of 3-Manifolds 
{\it Annals of Mathematics Studies, Princeton University Press} (1994).





\bibitem{KauffmanOgasasq}
L. H. Kauffman and  E. Ogasa: 
Steenrod square for virtual links toward Khovanov-Lipshitz-Sarkar
stable homotopy type for virtual links,  
   arXiv:2001.07789 [math.GT].



\bibitem{KOq}
L. H. Kauffman and E. Ogasa: 
 Quantum Invariants of Links and 3-Manifolds with Boundary defined via Virtual Links, 
arXiv 2108.13547 math.GT




\bibitem{KauffmanNikonovOgasa}
L. H. Kauffman, I. M. Nikonov, and E. Ogasa: 
Khovanov-Lipshitz-Sarkar homotopy type 
for links in thickened higher genus surfaces
arXiv: 2007.09241[math.GT].


\bibitem{KauffmanNikonovOgasaT2}
L. H. Kauffman, I. M. Nikonov, and E. Ogasa: 
Khovanov-Lipshitz-Sarkar homotopy type for links in thickened surfaces and those in $S^3$ with new modulis, 
arXiv:2109.09245  [math.GT]. 


\bibitem{Khovanov1999}
M. Khovanov.
\newblock {A categorification of the Jones polynomial}.
\newblock {\em Duke Mathematical Journal}, 101, 1999.





\bibitem{Kirbyc} 
R.  C. Kirby: A calculus for framed links in $S^3$
{\it Invent. Math.} 45 (1978), 35-56. 
 
 
\bibitem{Kirby}
R. C. Kirby:  The topology of 4-manifolds 
 {\it Lecture Notes in Math  (Springer Verlag) }
vol. 1374, 1989


\bibitem{KM}
R. C. Kirby and P. Melvin: 
The 3-manifold invariants of Witten and Reshetikhin-Turaev for sl(2, C)
{\it Inventiones Mathematicae} 105 (1991) 473–545. 



\bibitem{Kuperberg} 
G. Kuperberg: What is a virtual link? 
{\it Algebr. Geom. Topol.}  3 (2003) 587-591.




\bibitem{Lickorisho} 
W. B. R. Lickorish: 
 A representation of orientable combinatorial 3-manifolds 
{\it Ann. Math.} 76 (1962),
531-540.


\bibitem{Lickorish} 
W. B. R. Lickorish: 
Invariants for 3-manifolds from the combinatorics of the Jones polynomial 
{\it Pacific J. Math.} 149 (1991) 337-347.



\bibitem{Lickorishl} 
W. B. R. Lickorish: 
Three-manifolds and the Temperley-Lieb algebra 
{\it Mathematische Annalen} 290 (1991)   657–670. 




\bibitem{LSk} R. Lipshitz and S. Sarkar: A Khovanov stable homotopy type, 
{\it J. Amer. Math. Soc.} 27 (2014), no. 4, 983–1042. MR 3230817



\bibitem{LSs}
R. Lipshitz and S. Sarkar: 
A Steenrod square on Khovanov homology, 
{\it J. Topol.} 7 (2014), no. 3, 817–848. MR 3252965


\bibitem{LSr}
R. Lipshitz and S. Sarkar: 
A refinement of Rasmussen's s-invariant, 
{\it Duke Math. J.} 163 (2014), no. 5, 923–952. MR 3189434.



\bibitem{Man}
V O Manturov: Khovanov homology for virtual links with arbitrary coeficients, 
{\it Journal of Knot Theory and Its Ramifications} 16 (2007).



\bibitem{Mvbunrui} 
V. O. Manturov: 
Virtual Knots: The State of the Art,  
{\it $($in Russian$)$}
2010. 

\bibitem{MN}
V. O. Manturov and I. M. Nikonov: 
Homotopical Khovanov homology
{\it Journal of Knot Theory and Its Ramifications} 
24 (2015) 1541003. 







\bibitem{Igor}
I. M. Nikonov: 
Virtual index cocycles and invariants of virtual links, 
arXiv:2011.00248



\bibitem{Org} 
E. Ogasa: 
An elementary introduction to Khovanov-Lipshitz-Sarkar stable homotopy type, \\
https://www.researchgate.net/publication/352136009
$\_$
An
$\_$
elementary
$\_$
introduction
$\_$
to
$\_$
Khovanov-Lipshitz-Sarkar
$\_$
stable
$\_$
homotopy
$\_$
type \\
(The readers can find this article by typing in the title in a search engine.)


\bibitem{RT} 
N. Reshetikhin and V. G. Turaev: 
Invariants of 3-manifolds via link polynomials and quantum groups, 
{\it Inventiones mathematicae} 103 (1991) 547–597. 


 \bibitem{RT2}
 N. Reshetikhin and V. G. Turaev:
 Ribbon Graphs and their invariants derived from quantum groups,  
{\it Communications in Mathematical Physics},  127(1990) 1-26.
 


\bibitem{R}  
J. Roberts:  
Kirby calculus in manifolds with boundary 
{\it 
Proceedings of 5th Gökova. Geometry-Topology Conference.}21(1997)111-117, 
arXiv:math/9812086. 

\bibitem{Ru} W. Rushworth: Doubled Khovanov Homology, 
{\it Can. j. math.} 70 (2018) 1130-1172. 


\bibitem{Seed} C. Seed: 
Computations of the Lipshitz-Sarkar Steenrod square on Khovanov homology,
arXiv:1210.1882.




\bibitem{Tub}
D. Tubbenhauer: 
Virtual Khovanov homology using cobordisms, 
{\it J. Knot Theory Ramifications} 23 (2014), no. 9, 1450046, 91 pp. 



\bibitem{Viro} Viro; 
Khovanov homology of Signed diagrams 2006 (an unpublished note).





\bibitem{W} E. Witten:  Quantum field theory and the Jones polynomial
{\it Comm. Math. Phys.  } 121 (1989) 351-399. 
}

\end{thebibliography}
\end{document}